%% file: main.tex
\documentclass[10pt,twocolumn,twoside]{IEEEtran}

\usepackage{amsmath,amsfonts,amsthm,cite,times,enumitem,url,hyperref,bm, amssymb,color}
\newcommand{\tb}{\color{blue}}
\usepackage{bbold}
\usepackage{tcolorbox}

\usepackage{algorithm}
\usepackage[noend]{algorithmic}
\algsetup{indent=2em}

\usepackage{lipsum}
\usepackage{caption}
\usepackage{tikz,subcaption,cite}
\usetikzlibrary{decorations.pathmorphing}

\usepackage{enumitem}

\newtheorem{theorem}{Theorem}
\newtheorem{lemma}{Lemma}
\newtheorem{proposition}{Proposition}
\newtheorem{remark}{Remark}

\newtheorem{assumption}{Assumption}
\newtheorem{example}{Example}

\newcommand{\suchthat}{\;\ifnum\currentgrouptype=16 \middle\fi|\;}

\newcommand{\real}{\mathbb{R}}

\newcommand{\diag}{\mathrm{diag}}

\newcommand{\scirc}{\raise1pt\hbox{$\,\scriptstyle\circ\,$}}

\newcommand\oprocendsymbol{\hbox{$\square$}}
\newcommand\oprocend{\relax\ifmmode\else\unskip\hfill\fi\oprocendsymbol}

\definecolor{Gray}{gray}{0.9}

\DeclareSymbolFont{bbold}{U}{bbold}{m}{n}
\DeclareSymbolFontAlphabet{\mathbbold}{bbold}
\newcommand{\vect}[1]{\mathbbold{#1}}

\usepackage{lipsum}

\begin{document}
\title{\LARGE  \bf  
Network Critical Slowing Down: \\ Data-Driven Detection of Critical Transitions in Nonlinear Networks }

\author{Mohammad Pirani\IEEEauthorrefmark{1} and Saber Jafarpour\IEEEauthorrefmark{1}
\thanks{ M.~Pirani is with the Department of Mechanical and Mechatronics Engineering, University of Waterloo, E-mail: {\tt mpirani@uwaterloo.ca}. S. Jafarpour is with the Decision and Control Laboratory at Georgia Institute of Technology,  USA, E-mail: {\tt saber@gatech.edu}.}
\thanks{* These authors contributed equally.}%
}


\maketitle

\begin{abstract}
In a Nature article, Scheffer et al. presented a novel data-driven framework to predict critical transitions in complex systems. 
These transitions, which may stem from failures, degradation, or adversarial actions, have been attributed to bifurcations in the nonlinear dynamics.
Their approach was built upon the phenomenon of {\it critical slowing down}, i.e.,  slow recovery in response to small perturbations near bifurcations. 
We extend their approach to detect and localize critical transitions in  nonlinear networks. 
By introducing the notion of {\it network critical slowing down}, the objective of this paper is to detect that the network is undergoing a bifurcation only by analyzing its signatures from measurement data. 
We focus on two classes of widely-used nonlinear networks: (1) Kuramoto model for the synchronization of coupled oscillators, and (2)  attraction-repulsion dynamics in swarms, each of which presents a specific type of bifurcation. 
Based on the phenomenon of critical slowing  down, we study the asymptotic behavior of the perturbed system away and close to the bifurcation and leverage this fact to develop a deterministic method to detect and identify critical transitions in  nonlinear networks. 
Furthermore, we study the state covariance matrix subject to a stochastic noise process away and close to the bifurcation and use it to develop a stochastic framework for detecting critical transitions. 
Our simulation results show the strengths and limitations of the methods. 
\end{abstract}

\begin{IEEEkeywords}
Critical Slowing Down, Bifurcation, Nonlinear Coupled Oscillators, Attraction-Repulsion Dynamics. 
\end{IEEEkeywords}

\IEEEpeerreviewmaketitle

\section{Introduction}

\paragraph*{Motivation}

In complex systems, transitions from  stable to unstable operating modes due to the change in system's parameter are seen in a broad range of real-world applications, from biological and  ecological systems \cite{strogatz} to infrastructure systems \cite{Kwatny}. Examples include: (1) cluster synchronization in complex networks as a result of the infrastructure degradation over time \cite{Pecora}, (2) the onset of Alzheimer’s disease caused by changes in brain functional network connectivity affected by aging \cite{brown}, and (3) the onset of the epilepsy which is correlated to the gradual changes in the functional connectivity in the brain's temporal lobe \cite{Tracy}. In some application, these transitions plays crucial role in the operation of the network. However, in many other applications, they can endanger network's safety and potentially results in catastrophic events. From a mathematical perspective, these phenomena are described by bifurcations; changes in the stability of an equilibrium due to parameter variations. Detecting that a system is undergoing a bifurcation is not easy, because it usually requires knowledge about the system model and values of the bifurcating parameters.  However, it is shown that bifurcations usually leave traces on the system's measurements \cite{Scheffer}. Inspired by those observations,  we introduce data-driven algorithms to detect and localize critical transitions in two widely-studied nonlinear networks.

\paragraph*{Literature Review}
The  dynamics representing the evolution of natural systems, e.g., biological, ecological, and climate systems, represent different types of nonlinear phenomena, including the bifurcation. Thus, first attempts to use data in analyzing the behaviour of large-scale networks and detecting their critical transitions, belong to the researches in natural sciences and medicine \cite{Lenton, Litt}. One of the signatures of complex systems near bifurcation is the phenomenon of {\it critical slowing down}; a slow recovery in response to perturbations near a bifurcation. Critical slowing down was first introduced in statistical mechanics \cite{Fisher} and was  revisited in a variety of complex systems \cite{Leemput, Lenton}. The effect of critical slowing down  on the autocorrelation of the system's state, its variance, and higher order moments of data distribution was discussed in \cite{Dakos, Scheffer} which became useful in detecting tipping points in complex systems using data. However, these approaches are usually based on the analysis of an aggregated scalar dynamical system and do not consider the role of network structure in these critical transitions. In a separate line of research, detection and identification of abrupt changes in dynamical systems has been extensively studied in engineering systems, see \cite{Basseville} and references therein. These methods mainly consist of two steps: (1) generating a residual, i.e., a signal showing the deviation of the measurement from the system's output in the normal condition, and (2) designing a decision rule based upon these residuals. The above steps require an extensive knowledge about the system's dynamics in the normal and faulty conditions and are usually applicable to linear systems.


In nonlinear networks, not only the bifurcation can cause instability in the system, the underlying network topology can facilitate the propagation of this failure throughout the  system. This phenomenon is known as  cascading failure that can potentially lead to the collapse of the entire infrastructure \cite{Dobson, RMD:17}. Hence, early detection and localization of failures in a network can help preventing their dissemination to the healthy parts. 
In this paper, we study two well-known dynamical systems, as discussed below:

\paragraph*{\textbf{Synchronization of Coupled Oscillators}}
Synchronization is a pervasive phenomenon in many engineering and natural systems including biological systems~\cite{TV-AC-EB-IC-OS:95}, robotics~\cite{DAP-NEL-RS-DG-JKP:07}, brain networks~\cite{PAT:03}, circuits and radio technology \cite{VanDerPol}, and power grids~\cite{florian}. 
One of the simplest model to study synchronization in oscillator networks is Kuramoto oscillator model. 
It is well-known that Kuramoto oscillators can undergo a bifurcation from synchrony to incoherency.  Indeed, for Kuramoto oscillator networks, synchronization is determined by a trade-off between coupling strength and oscillators' heterogeneity~\cite{florian}. A large body of work in the literature is devoted to estimating or characterizing this trade-off~\cite{jadbabaie}~\cite{florian}~\cite{saber}. 
In many real world applications of oscillator networks, finding the exact onset of bifurcation is a critical task for operation and safety of these networks. One of the most visible examples is the power grids which, as a result of the unprecedented penetration of renewable energy resources and large increase in power demand, are being pushed toward their maximum capacity. In this context, understanding the threshold of bifurcation can help system operators to ensure safety and security of the frequency synchronization. Another example is the brain network where computing the onset of bifurcation is important to understand different modes of functionality of the brain.


\paragraph*{\textbf{Attraction-Repulsion Dynamics}} The attraction-repulsion functions are widely used in formation and swarming of robots \cite{gazi, spears, Schwager}. They are inspired from animal grouping and biological systems with attraction having longer range than repulsion \cite{Okubo, Lazarus}. The dynamics of these swarms are governed by an interplay between these two forces and the network reaches local stability \cite{gazii}.  In these systems, a bifurcation can happen when the repulsion dominates the attraction, resulting in a change in the shape of the network.  An example  is the cluster formation of particles in granular flows where the cluster size is the resultant of the two forces \cite{zohdi}. Another example is the tissue growth in multi-cellular organisms in which density-dependent inhibition acts as a repulsion force between cells in a tissue and inhibits them to grow in a bounded environment \cite{stoker}. This density dependent inhibition is usually applied by cell-cell contacts. Tumor cells have often lose this density-dependent inhibition, resulting in a non-stopping cell production. Hence, understanding the logic behind bifurcations in those systems can help us to detect and identify them faster. 

\paragraph*{Contributions}
In this paper, we develop two data-driven methods (a deterministic and a stochastic) to detect critical transitions in nonlinear network systems. The contributions of this paper are as follows: 
\begin{itemize}
    \item We study bifurcation for coupled oscillators and  attraction-repulsion dynamics with two agents. Using suitable reduced-order models, we characterize the behavior of these systems at  bifurcation and their bifurcation type. Following~\cite{Scheffer}, we observe the phenomenon of critical slowing down near bifurcations and applied deterministic and stochastic methods to detect the onset of bifurcation.  


    \item For networks of coupled oscillators and networks of attraction-repulsion dynamics, we investigate asymptotic behaviors of parameterized nonlinear networks and provide sufficient conditions, in terms of network parameters and structure, under which a bifurcation happens (Propositions~\ref{prop:saddle-node} and \ref{prop:attraction}).

    \item  For networks of coupled oscillators and networks of attraction-repulsion dynamics, we study the signature of the bifurcation in these dynamics using a deterministic and a stochastic approach. In the deterministic approach, we investigate the asymptotic behavior of nonlinear networks away and close to the bifurcation (Theorem~\ref{thm:asymbehavior-kuramoto}). Building upon this result, we introduce the notion of network critical slowing down and develop a completely data-driven heuristic algorithm to detect bifurcations in the network and localize the bifurcating edges. We further demonstrate a dichotomy between the asymptotic behaviour of the state covariance before and at the bifurcation (Theorem~\ref{thm:nonlinearoscstoch}) and use it as a stochastic indicator of the bifurcation in the network. 
    
    \item We demonstrate the efficiency of our algorithms using several numerical simulations and compare the abilities, limitations, and application domains of the proposed deterministic and stochastic methods.

    
\end{itemize}

Finally, we emphasize that the objective of this work is to use data-driven methods to directly detect bifurcation without any need to identify system's parameters. We do not use system identification techniques and, thus, prevailing conditions to identify the system's parameters (i.e., persistence of excitation) do not need to be satisfied.

\section{Modeling and Problem Statement}
\label{sec:formulaaa}
In this section, we introduce two classes of nonlinear dynamics, namely attraction-repulsion dynamics and synchronization of coupled oscillators, and then state the problem.

\subsection{Notation} \label{sec:definitions}

Given a vector $v\in \real^n$, we define the diagonal matrix $[v]\in \real^{n\times n}$ by $[v]_{ii} = v_i$, for every $i\in \{1,\ldots,n\}$. For a symmetric matrix $M$, the eigenvalues are ordered as  $\lambda_1(M) \le \lambda_2(M) \le \ldots \le \lambda_n(M)$. We denote an undirected weighted graph by  $\mathcal{G}=\{\mathcal{V},\mathcal{E}, \mathcal{W}\}$,  where $\mathcal{V}\in \{1,2,..., n\}$ is a set of nodes  and $\mathcal{E} \subset \mathcal{V}\times\mathcal{V}, |\mathcal{E}|=m$ is the set of $m$ edges, and $\mathcal{W}= \{w_{ij}\}_{(i,j)\in \mathcal{E}}$ represents the edge weights.  The neighbors of node $i \in \mathcal{V}$ are given by the set $\mathcal{N}_i = \{j \in \mathcal{V} \mid (i, j) \in \mathcal{E}\}$. The  adjacency matrix of the graph is a symmetric and binary $n \times n$  matrix $A$, where element $A_{ij}=w_{ij}\in \mathcal{W}$ if $(i, j) \in \mathcal{E}$ and zero otherwise. The degree of node $i$ is $d_i=\sum_{j=1}^nA_{ij}$.  The Laplacian matrix of the graph is $L \triangleq D - A$, where $D = \diag(d_1, d_2, \ldots, d_n)$.  Given an undirected weighted graph $\mathcal{G}=\{\mathcal{V},\mathcal{E}, \mathcal{W}\}$, the eigenvalues of its Laplacian are real and nonnegative, and are denoted by $0 = \lambda_1(\mathcal{G}) \le \lambda_2(\mathcal{G}) \le \ldots \le \lambda_n(\mathcal{G})$.  An orientation of the graph $\mathcal{G}$ is defined by assigning  a direction (arbitrarily) to each edge in $\mathcal{E}$. For graph $\mathcal{G}$ with $m$ edges, numbered as $e_1, e_2, ..., e_m$, its node-edge incidence matrix $B(\mathcal{G})\in \mathbb{R}^{n\times m}$, or simply $B$, is defined as 
$$[B]_{kl}=
  \begin{cases}
    1       & \quad  \text{if node $k$ is the head of edge $l$},\\
   -1  & \quad \text{if node $k$ is the tail of edge $l$},\\
   0  & \quad \text{otherwise}.\\
  \end{cases}
  $$
The graph Laplacian satisfies $L=BWB^{\top}$, where $W$ is a diagonal matrix representing the edge weights.
Given a set of nodes $S\subseteq \mathcal{V}$, the cutset of $S$ is defined as $\partial S=\{(u,v)\in\mathcal{E}|u\in S, v\in \mathcal{V}\setminus S\}$. A $2$-cutset is a cutset whose removal splits the graph into two connected components.\footnote{ In this paper, for the sake of brevity, we have developed our results considering 2-cutsets. However, our results can be easily generalized to general cutsets.} Given a set $S\subseteq \mathcal{V}$, its indicator vector $\chi^S\in \real^n$ is defined by:
\begin{align*}
    \chi^S_i =\begin{cases}
      +1 & i\in S,\\
      -1 & i\ne S.
    \end{cases}
\end{align*}
Let $\mathcal{G}(\alpha)=\{\mathcal{V},\mathcal{E},\mathcal{W}(\alpha)\}$ be a family of undirected weighted graphs parameterized by $\alpha\in \mathcal{I}$. We define the associated \emph{lower-bound graph} $\underline{\mathcal{G}} = \{\mathcal{V},\mathcal{E},\underline{\mathcal{W}}\}$ where $\underline{\mathcal{W}} = \{\underline{w}_{ij}\}_{(i,j)\in \mathcal{E}}$ is defined by $\underline{w}_{ij} = \inf_{\alpha\in \mathcal{I}} w_{ij}(\alpha)$ for every $(i,j)\in \mathcal{E}$.


\subsection{Nonlinear Coupled Oscillators}

 One of the simplest model for studying network of nonlinear coupled oscillators is the Kuramoto model. In the Kuramoto model, the state of each agent is its phase $\theta_i\in \mathbb{R}$ and evolves according to the following first-order dynamics
\begin{equation}\label{eq:kuramoto}
    \dot{\theta}_i=\omega_i(\alpha)-\sum_{j\in \mathcal{N}_i}^n a_{ij}\sin (\theta_i-\theta_j), \quad i=1,2,...,n
\end{equation}
where $a_{ij}$ is the coupling strength between agents $i$ and $j$ and $\omega_i(\alpha)$ is the natural frequency of agent $i$, parameterized by $\alpha\in (0,\infty)$. Recalling $|\mathcal{E}|=m$, let $\mathcal{A}\in \real^{m\times m}$ be a diagonal edge weight matrix in which $\mathcal{A}_{ee} = a_{ij}$ where $e=(i,j)$ and $B\in \real^{n\times m}$ be the incidence matrix of $\mathcal{G}$. Then the coupled oscillator dynamics~\eqref{eq:kuramoto} can be written in the matrix form:
\begin{align}\label{eq:coupled-oscillator}
    \dot{\theta} = \omega(\alpha) - B\mathcal{A}\sin(B^{\top}\theta),
\end{align}
where $\theta=(\theta_1,\ldots,\theta_n)^{\top}\in \real^n$ and $\omega(\alpha)=(\omega_1(\alpha),\ldots,\omega_n(\alpha))^{\top}\in \real^n$ are the vector of states and natural frequencies respectively. 


\subsection{Attraction-Repulsion Dynamics in Swarms}
A swarm is a network of $n$ agents in a $k$-dimensional Euclidean space. The state of each agent is denoted by $x_i\in \mathbb{R}^k$ and evolves according to the following dynamics 
\begin{equation}\label{eqn:n08g86r}
    \dot{x}_i=\sum_{(i,j)\in \mathcal{E}}^n g(x_i-x_j), \quad i=1,2,...,n
\end{equation}
where $g:\mathbb{R}^k\to\mathbb{R}^k$ is a nonlinear function modeling the attraction and repulsion forces between the individuals. Among several classes of attraction-repulsion functions $g(\cdot)$,  linear attraction and exponential repulsion is widely used in the literature \cite{gazi, Brambilla}, inspired from natural swarms with long-range attraction  and short-range repulsion. Using this model, the dynamics of each agent become
\begin{equation}\label{eqn:1asf}
  \dot{x}_i=-\sum_{(i,j)\in \mathcal{E}}(x_i-x_j)\Big(w_a^{ij}(\alpha)-w_r^{ij}\exp\big(\tfrac{-\|x_i-x_j\|^2}{c}\big)  \Big),
\end{equation}
where $w_a^{ij}$ is the attraction coefficient and $w_r^{ij}$ is the repulsion coefficient and they represent the strength of attraction and repulsion forces between the node pair $(i,j)$, respectively.
Here, $w_a^{ij}$ is parameterized by the parameter $\alpha\in (0,\infty)$. 
One can write~\eqref{eqn:1asf} in the following matrix form: 
\begin{equation}\label{eqn:nougf86}
    \dot{x}=-B\overline{\mathcal{A}}(x,\alpha)B^{\top}x
\end{equation}
where $B$ is the incidence matrix of the graph and $\overline{\mathcal{A}}=\diag \left(w_a^{ij}(\alpha)-w_r^{ij} \exp (\frac{-\|x_i-x_j\|^2}{c}) \right)$ for every $(i,j)\in \mathcal{E}$. 


\subsection{Problem Statement}
Consider a dynamical system on the graph $\mathcal{G}$ with the parameterized vector field $f:\real^n\times \real_{\ge 0}\to \real^n$, in which the states of nodes evolve as follows: 
\begin{equation}\label{eqn:mainproblem}
\dot{x}=f(x,\alpha),
\end{equation}
where   $x\in \real^n$ is the state of system, $\alpha\in (0,\infty)$ is the bifurcation parameter, and the parameterized vector field $f$ take either form of \eqref{eq:coupled-oscillator} or \eqref{eqn:nougf86}. We also assume that:
\begin{itemize}
    \item We have access to state measurements $x$.
    \item We have no knowledge of the network parameters and the bifurcation parameter $\alpha$. 
\end{itemize}
 We know that parameter $\alpha$ is changing and a bifurcation happens at some parameter value $\alpha=\alpha^*$. The problem is to detect when the system is close to bifurcation, i.e.,  $\alpha$ is close to $\alpha^*$. We address the above problem by studying the response of the system \eqref{eqn:mainproblem} to perturbations. The perturbations can be considered to be deterministic (i.e., disturbance signals with known magnitude) or stochastic (i.e., Gaussian noise with a known mean and covariance).

\section{Bifurcation Analysis for 2-dimensional Systems}

For the sake of simplicity in the exposition, in this section and Section \ref{sec:detection-2}, we start with analyzing the bifurcation in networks with two agents. Our two-dimensional analysis can be also useful when studying aggregation of agents in a network or model-order reduction \cite{Chow}.

\subsection{Nonlinear Coupled Oscillators}
\label{sec:cngt6ruhv29}
\begin{figure*}[t!] 
\centering
\includegraphics[scale=.60]{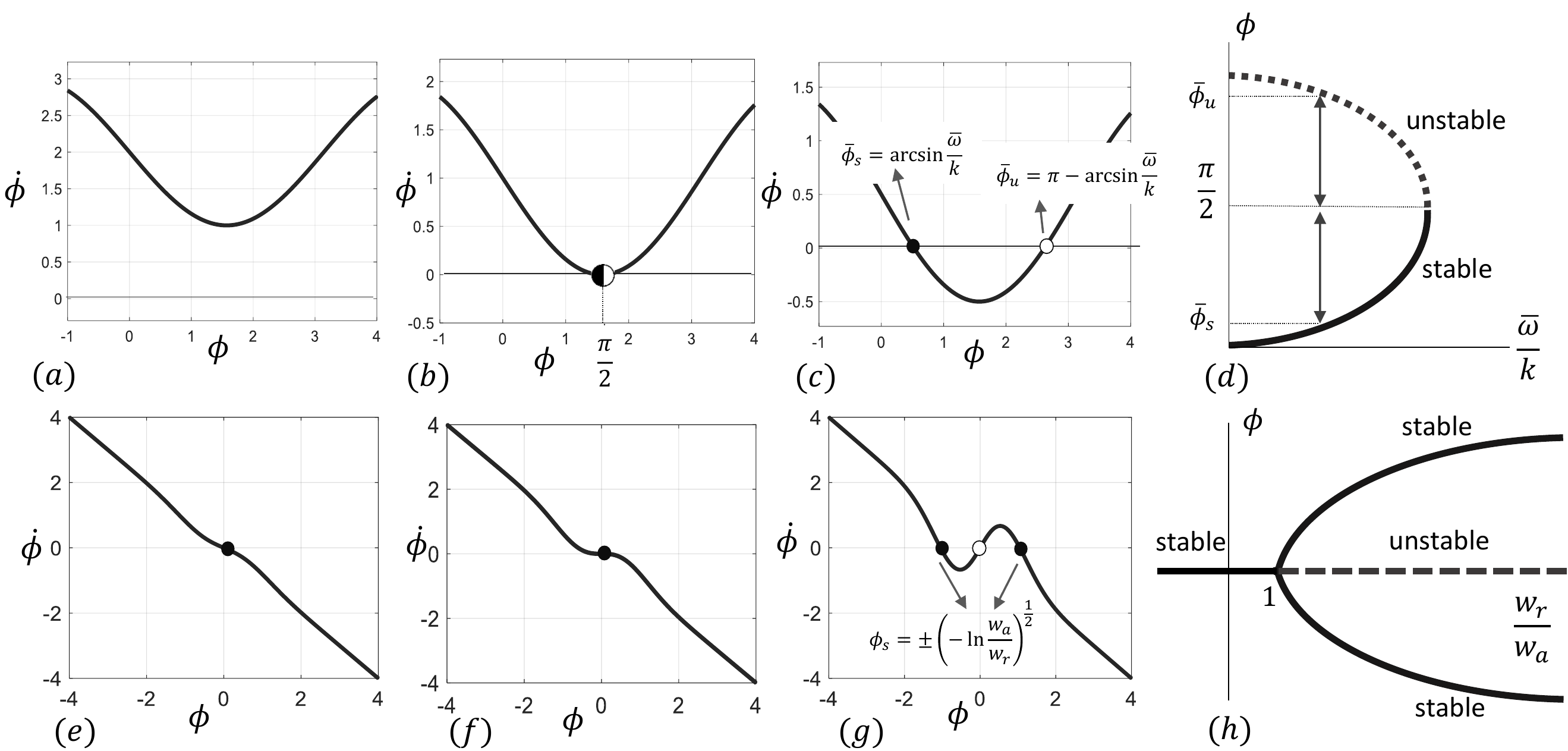}
\caption{\footnotesize{Top: Phase plane for reduced-order coupled oscillator dynamics~\eqref{eq:reduced-order-coupled-oscillator} with different values of $\frac{\bar{\omega}}{k}$. (a) $\frac{\bar{\omega}}{k}=2$ (b) $\frac{\bar{\omega}}{k}=1$ (c) $\frac{\bar{\omega}}{k}=0.5$ (d) Bifurcation diagram. Bottom: Phase plane for reduced-order attraction repulsion dynamics~\eqref{eqn:cn08hgt} with different values of $\frac{w_r}{w_a}$ and $c=1$. (e) $\frac{w_r}{w_a}=0.5$, (f) $\frac{w_r}{w_a}=1$ (g) $\frac{w_r}{w_a}=3$ (h) Bifurcation diagram.}}
\label{fig:sst968svgvnj}
\end{figure*}

Consider a network consists of two oscillators
\begin{align}\label{eqn:cman09hg7}
\dot{\theta}_1=\omega_1(\alpha)-k\sin(\theta_1-\theta_2),\nonumber\\
\dot{\theta}_2=\omega_2(\alpha)-k\sin(\theta_2-\theta_1).
\end{align}
where $\omega_i(\alpha)$ is the natural frequency of node $i$, parameterized by the bifurcation parameter $\alpha\in \real$ and $k$ is a fixed coupling gain. We define $\phi=\theta_1-\theta_2$ and $\bar{\omega}(\alpha)=\omega_1(\alpha)-\omega_2(\alpha)$. We show that when the ratio between the difference in natural frequencies and the coupling, i.e., $\frac{\overline{\omega}(\alpha)}{k}$, exceeds a threshold, a \emph{saddle node bifurcation} happens.\footnote{See \cite{strogatz} for details about several types of bifurcations in nonlinear systems.} To show this, we subtract the two equations in \eqref{eqn:cman09hg7} to get the \emph{ reduced-order coupled oscillator dynamics} as:
\begin{align}\label{eq:reduced-order-coupled-oscillator}
\dot{\phi}=\bar{\omega} (\alpha) -k\sin\phi.
\end{align}
 For $\frac{\bar{\omega}(\alpha)}{k}>1$, the reduced-order dynamics does not have any equilibrium point, as shown in Fig. \ref{fig:sst968svgvnj} (a).  When $\frac{\bar{\omega}(\alpha)}{k}=1$, a saddle node bifurcation  happens at $\phi=\frac{\pi}{2}$  (Fig. \ref{fig:sst968svgvnj} (b)). For $\frac{\bar{\omega}(\alpha)}{k}<1$, there are two equilibria, one stable and one unstable (Fig. \ref{fig:sst968svgvnj} (c)).\footnote{As an example, in power systems the natural frequencies of the generators change over time, due to degradation or power demand~\cite{canizares}. } 

\subsection{Attraction-Repulsion Dynamics}\label{sec:mn09hb82b}

Consider two agents with scalar states $x_1$ and $x_2$. The attraction-repulsion dynamics are given by
\begin{align}\label{eq:attract-repulsion}
     \dot{x}_1=-(x_1-x_2)\bigg(w_a(\alpha)-w_r.\exp\big(\frac{-\|x_1-x_2\|^2}{c}\big)  \bigg),\nonumber\\
     \dot{x}_2=-(x_2-x_1)\bigg(w_a(\alpha)-w_r.\exp\big(\frac{-\|x_2-x_1\|^2}{c}\big)  \bigg).
\end{align}
Subtracting the two equations and defining $\phi=x_1-x_2$, gives the \emph{reduced-order attraction-repulsion dynamics} as follows:
\begin{align}\label{eqn:cn08hgt}
     \dot{\phi}=-2\phi\bigg(w_a(\alpha)-w_r.\exp\big(\frac{-\|\phi\|^2}{c}\big)  \bigg).
\end{align}
 It can be shown that when the ratio between the repulsion coefficient and the attraction coefficient exceeds a threshold, a {\it supercritical pitchfork bifurcation} happens and the system transitions from one fixed point to three fixed points.  For $\frac{w_r}{w_a(\alpha)}\leq 1$, the origin is a stable equilibrium point. For  $\frac{w_r}{w_a(\alpha)}> 1$, the origin becomes unstable with two symmetrical stable equilibria at $\phi_e=\pm\big(-c\ln(\frac{w_a}{w_r}) \big)^{\frac{1}{2}}$, as shown in Fig.~\ref{fig:sst968svgvnj} (e-g). In this figure, stable and unstable equilibrium points are shown with black and white circles, respectively.  Figure \ref{fig:sst968svgvnj} (h) shows the bifurcation diagram of \eqref{eqn:cn08hgt}.

\section{Detecting Critical Transitions for 2-dimensional systems}\label{sec:detection-2}

We focus on two different approaches to detect bifurcations in nonlinear systems: i) deterministic approach by considering the recovery rate after perturbations, and ii) stochastic approach by considering the autocorrelation and variance of the state. Our framework is built upon the phenomenon of critical slowing down. It refers to the tendency of a system to take longer to return to equilibrium after perturbations as discussed in~\cite{Scheffer}. 

\subsection{Deterministic Method}\label{subsec:recovery-2}

 In this method, we consider the recovery rate of the system subject to a known and deterministic perturbation to detect bifurcations. For a given parameter $\alpha$, consider the dynamical system~\eqref{eqn:mainproblem} and suppose that 
 $\overline{x}_{\alpha}$ is a locally asymptotically stable equilibrium point of~\eqref{eqn:mainproblem}. We study the behavior of ~\eqref{eqn:mainproblem} around $\overline{x}_{\alpha}$ under a small perturbation $\epsilon$ and get: 
\begin{equation}\label{eqn:vm089y}
    \frac{d(\overline{x}_{\alpha}+\epsilon)}{dt}=f(\overline{x}_{\alpha}+\epsilon).
\end{equation}
Linearizing  \eqref{eqn:vm089y} using a first-order Taylor expansion yields 
\begin{equation}\label{eqn:vm015say}
    \frac{d(\overline{x}_{\alpha}+\epsilon)}{dt}=f(\overline{x}_{\alpha}+\epsilon)= f(\overline{x}_{\alpha})+\frac{\partial f}{\partial x}\bigg|_{\overline{x}_{\alpha}}\epsilon + \mathcal{O}(\epsilon^2).
\end{equation}
 Therefore, the dynamics for the perturbation $\epsilon$ around the equilibrium point $\overline{x}_{\alpha}$ can be approximated by the associated perturbation dynamics
\begin{equation}\label{eqn:perturbation_dynamics}
    \frac{d\epsilon}{dt}=\frac{\partial f}{\partial x}\bigg|_{\overline{x}_{\alpha}}\epsilon.
\end{equation}
 By~\cite[Theorem 4.7]{khalil}, if the origin is the exponentially stable equilibrium point of~\eqref{eqn:perturbation_dynamics}, then  $\overline{x}_{\alpha}$ is an locally exponentially stable equilibrium point of~\eqref{eqn:vm089y}. 
 Now we study the perturbation dynamics~\eqref{eqn:perturbation_dynamics} for the reduced order models of coupled oscillators dynamics~\eqref{eq:reduced-order-coupled-oscillator} and attraction-repulsion dynamics~\eqref{eqn:cn08hgt} and propose algorithms for detection of critical transitions in these systems.

\subsubsection*{Nonlinear Coupled Oscillators}

For the reduced-order coupled oscillator~\eqref{eq:reduced-order-coupled-oscillator}, bifurcation occurs when $\frac{\bar{\omega}(\alpha)}{k}=1$. For $\frac{\bar{\omega}(\alpha)}{k}<1$, the dynamical system~\eqref{eq:reduced-order-coupled-oscillator} has two equilibrium points  $\overline{\phi}_s=\mathrm{arcsin}(\frac{\bar{\omega}(\alpha)}{k})$ and $\overline{\phi}_u=\pi-\mathrm{arcsin}(\frac{\bar{\omega}(\alpha)}{k})$ as shown in Fig.~\ref{fig:sst968svgvnj} (c). It is easy to see that $\overline{\phi}_s=\mathrm{arcsin}(\frac{\bar{\omega}(\alpha)}{k})$ is the locally asymptotically stable equilibrium point of~\eqref{eq:reduced-order-coupled-oscillator}. Using \eqref{eqn:perturbation_dynamics}, the perturbation dynamics around $\overline{\phi}_s$ can be approximated by:
\begin{equation}\label{eqn:vm0s648egay}
    \frac{d\epsilon}{dt}= -\cos(\mathrm{arcsin}(\tfrac{\bar{\omega}(\alpha)}{k}))\epsilon.
\end{equation}
For $\frac{\bar{\omega}(\alpha)}{k}<1$, $\epsilon=0$ is an exponentially stable equilibrium point of \eqref{eqn:vm0s648egay} and thus the perturbation $\epsilon$ converges to zero. 
%
According to~\eqref{eqn:vm0s648egay}, by approaching to the bifurcation, i.e, $\frac{\bar{\omega}(\alpha)}{k}\to 1$, the rate of the recovery of perturbation $\epsilon$ diminishes. 
%

\subsubsection*{Attraction-Repulsion Dynamics} For the reduced-order attraction-repulsion dynamics~\eqref{eqn:cn08hgt} bifurcation occurs when $\frac{w_r}{w_a(\alpha)}=1$. For $\frac{w_r}{w_a(\alpha)}\le 1$, the only equilibrium point of~\eqref{eqn:cn08hgt} is $\overline{\phi}=0$ and thus \eqref{eqn:perturbation_dynamics} becomes 
\begin{equation}\label{eqn:vm1658wgy}
    \frac{d\epsilon}{dt}= 2(w_r-w_a(\alpha))\epsilon.
\end{equation}
For $\frac{w_r}{w_a(\alpha)}< 1$, based on \eqref{eqn:vm1658wgy}, the perturbation goes exponentially to zero, i.e., $\epsilon(t)=\epsilon (0)e^{2(w_r-w_a(\alpha))t}$. 
By approaching to the bifurcation, i.e., $\frac{w_r}{w_a(\alpha)}\to 1$, the rate of the recovery after perturbation $\epsilon$ diminishes. 
%

\begin{example}\label{exp:11}
Consider nonlinear coupled oscillator model~\eqref{eq:reduced-order-coupled-oscillator} with the coupling $k=2$. A constant perturbation signal $\epsilon=0.4$ is applied to state $\phi$ of the system from $t=6s$ to $t=8s$.  The perturbation signal before bifurcation for different values of $\frac{\bar{\omega}}{k}$ is shown in Fig.~\ref{fig:2majwew}, left. As shown in this figure, the closer value of $\frac{\bar{\omega}}{k}$ to the bifurcation point, i.e., $\frac{\bar{\omega}}{k}=1$, the larger the time required to return the perturbation to zero. The response at the bifurcation, $\frac{\bar{\omega}}{k}=1$, is shown in Fig.~\ref{fig:2majwew}, right which shows that the perturbation is no longer recovered.
 \begin{figure}[t!]
\centering
\includegraphics[scale=0.33]{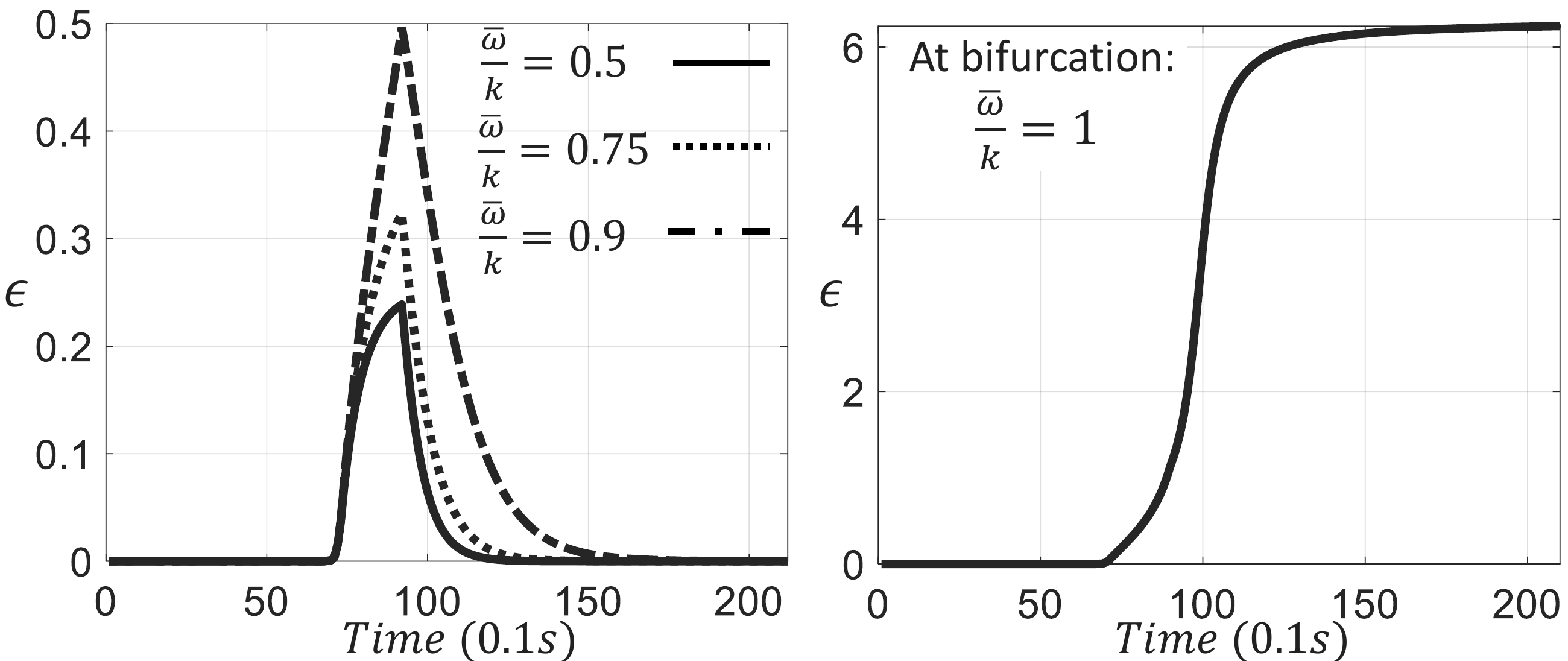}
\caption{\footnotesize{Nonlinear coupled oscillators: The evolution of $\epsilon$ before bifurcation, left, and during bifurcation, right.}}
\label{fig:2majwew}
\end{figure}

Now, consider attraction-repulsion dynamics~\eqref{eqn:cn08hgt}. The evolution of the perturbation signal $\epsilon$ for the attraction-repulsion dynamics is shown in Fig.~\ref{fig:2m1nn6} for different values of $\frac{w_r}{w_a}$ and  $c=1$. 
\begin{figure}[t!]
\centering
\includegraphics[scale=0.31]{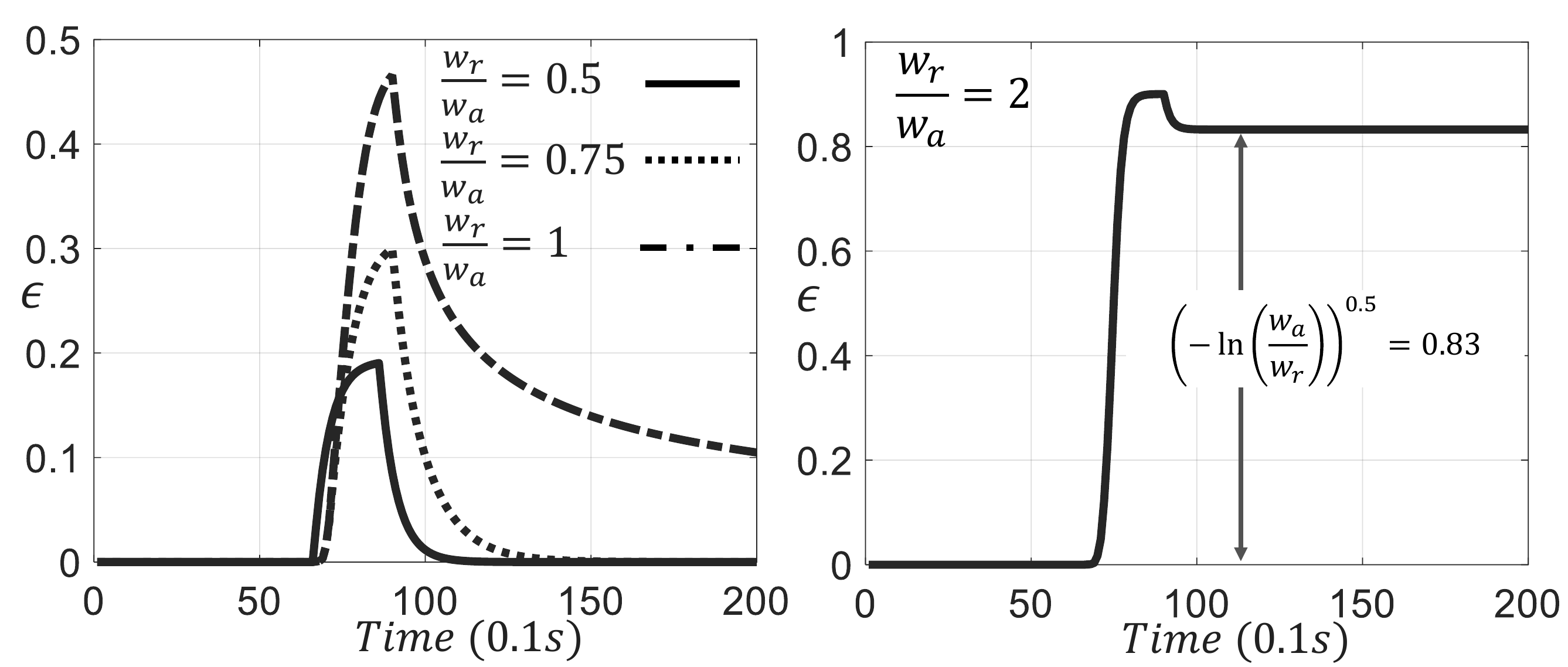}
\caption{\footnotesize{Attraction-repulsion dynamics: The evolution of $\epsilon$ before bifurcation, left, and after bifurcation, right.}}
\label{fig:2m1nn6}
\end{figure}
The bifurcation occurs at $\frac{w_r}{w_a} =1$. Similar to the coupled oscillator dynamics, by approaching to the bifurcation, the tendency of the system to recover the perturbation decreases. However, during and after bifurcation,  the attraction-repulsion dynamics present a different behaviour. In the bifurcation, the origin is still a globally stable equilibrium (unlike saddle-point equilibrium in the coupled oscillator dynamics). Thus, as shown in Fig.~\ref{fig:2m1nn6}, left, the perturbation will eventually recovered to the origin. After bifurcation, Fig.~\ref{fig:2m1nn6}, right, the origin becomes unstable and the perturbation is recovered to one of the emerging stable equilibria, i.e., those shown in Fig.~\ref{fig:sst968svgvnj} (g).
\end{example}

\subsection{Stochastic Method} \label{sec:variancee}

In this section, we look at the autocorrelation of system's state subject to a stochastic noise to detect bifurcations. The intuition is that near the bifurcation, the rates of change of the system's states decrease. As a result, the state of the system at any given moment becomes more similar to its past state which leads to an increase in the  autocorollation of system's states.
To formally show this, we approximate the dynamics of the system close to the equilibrium $\overline{x}_{\alpha}$ by
\begin{align}\label{eqn:buf97g8}
    \frac{dx}{dt} = f(x) = f(\overline{x}_{\alpha}) +\frac{\partial f}{\partial x}\bigg|_{\overline{x}_{\alpha}} (x-\overline{x}_{\alpha}) + \mathcal{O}(\epsilon^2). 
\end{align}
Assume that a stochastic disturbance is applied in specific time steps $\Delta t$, for $\Delta t\in \mathbb{Z}_{\ge 0}$ to the above dynamics.
The state of the system after applying the stochastic white noise $\{\xi_{t}\}_{t=1}^{\infty}$ with variance $\sigma$ is given by 
 \begin{equation}\label{eqn:cm09h9o8g}
x_{t+1}-\overline{x}_{\alpha}=e^{\lambda \Delta t}(x_{t}-\overline{x}_{\alpha})+\xi_t,
 \end{equation}
Following the argument in Section~\eqref{subsec:recovery-2}, for  systems~\eqref{eq:reduced-order-coupled-oscillator} and~\eqref{eqn:cn08hgt}, the (approximated) exponential rates of recovery are 
\begin{align}
    \lambda&=-\cos(\mathrm{arcsin}(\tfrac{\bar{\omega}(\alpha)}{k})) && \text{Nonlinear Coupled Oscillator} \nonumber\\
    \lambda&=2(w_r-w_a(\alpha)) && \text{Attraction-Repulsion} 
\end{align}
 We define $e_t=x_t-\overline{x}_{\alpha}$. If $\lambda$ and $\Delta t$ are independent of $x_t$, one can interpret \eqref{eqn:cm09h9o8g} as a lag-1 autoregressive process, 
\begin{equation}\label{eqn:nc0a89v8}
    e_{t+1}=\gamma e_t + \xi_t.
\end{equation}
If the autocorrelation parameter $\gamma=e^{\lambda\Delta t}$ is close to zero, the state $x_{t}$ inherits the white noise characteristics of $\xi_t$ and if it is close to one, $e_{t+1}$ becomes closer to a red (autocorrelated) noise~\cite{Florescu}. Since $|\gamma|<1$, for stable system, the mean $\mathbb{E}(e_t)$ is identical for all values of $t$ by the definition of wide sense stationarity and the variance is calculated as \cite{Florescu}
\begin{equation}
   \lim_{t\to \infty} {\rm var}(e_t)=\frac{\sigma^2}{1-\gamma^2}.
\end{equation}
Close to the bifurcation, the recovery rate to the equilibrium is small, implying that $\lambda$ approaches zero. Thus, the autocorrelation $\gamma$ tends to one and the variance tends to infinity. 

\begin{figure}[t!]
\centering
\includegraphics[scale=0.35]{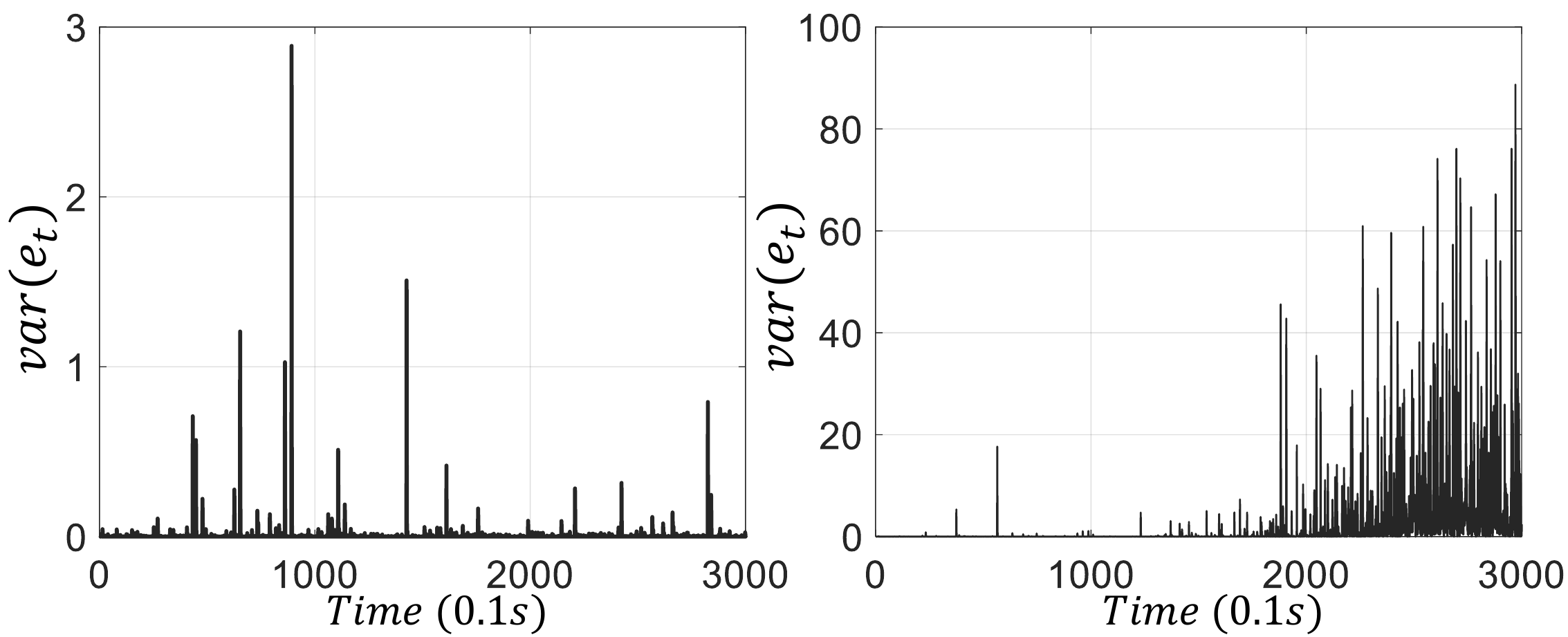}
\caption{\footnotesize{State variance away from (left) and close to the bifurcation (right).}}
\label{fig:2mbkvik6}
\end{figure}

\begin{example}\label{exp:12}
Consider nonlinear coupled oscillator~\eqref{eq:reduced-order-coupled-oscillator} with the coupling $k=2$. A zero-mean white noise with  variance $\sigma=5$ is applied to the system at every integer instance of time $t\in \mathbb{Z}_{\ge 0}$. The moving variance of the state is shown in Figure~\ref{fig:2mbkvik6}. 
\end{example}

So far, we focused on the critical transitions in scalar reduced-order models, ignoring the role of network structure. In the coming sections, we extend this framework to detect and localize bifurcations in nonlinear networks.

\section{Bifurcation Analysis in Networks}
\label{sec:network}

Before extending the detection methods to nonlinear networks, we first characterize the attracting sets of the each nonlinear network before and after the bifurcation.

\subsection{Nonlinear Coupled Oscillators}

Consider the coupled oscillator dynamics~\eqref{eq:coupled-oscillator} over a network $\mathcal{G}=\{\mathcal{V},\mathcal{E}, \mathcal{W}\}$. We suppose the vector of natural frequencies $\omega(\alpha)=[\omega_1(\alpha) \hspace{1mm} \omega_2(\alpha) \hspace{1mm}...\hspace{1mm}\omega_n(\alpha)]^{\top}$ is parameterized by a real-valued parameter $\alpha\in [0, \infty )$ satisfying the following assumption.

\begin{assumption}\label{ass:bifurcation}
There exist a $2$-cutset $\partial S$ for some $S\subseteq \mathcal{V}$ and an $\alpha^*\in (0, \infty)$ such that 
\begin{enumerate}
    \item for every $\alpha \in [0,\alpha^*)$ we have $\|B^{\top}L^{\dagger}\omega(\alpha)\|_{\infty}<1$;
    \item if $\alpha=\alpha^*$, then  $\left|(B^{\top}L^{\dagger}\omega(\alpha^*))_e\right|= 1$, for every $e\in \partial S$;
    \item for every $\alpha\in (\alpha^*,\infty)$,
    \begin{align*}
        \left|(B^{\top}L^{\dagger}\omega(\alpha))_e\right|&> 1, \mbox{ for every }  e\in \partial S, \\
        \left|(B^{\top}L^{\dagger}\omega(\alpha))_e\right|&< 1, \mbox{ for every } e\in \mathcal{E}\setminus\partial S.
    \end{align*}
\end{enumerate}
\end{assumption}

The term $B^{\top}L^{\dagger}\omega(\alpha)\in \real^m$ has a physical interpretation in flow networks. For acyclic undirected graphs, any single edge $e$ is an edge cut and the vector  $B^{\top}L^{\dagger}\omega(\alpha)\in \real^m$ is interpreted as the normalized flow through the edges of the network\cite{saber}. As a result, for acyclic graphs, Assumption \ref{ass:bifurcation} indicates that, at $\alpha=\alpha^*$, the flows on the $2$-cutset $\partial S$ reach their maximum value. The next proposition shows that the nonlinear coupled oscillator~\ref{eq:coupled-oscillator} undergoes a bifurcation at $\alpha=\alpha^*$ and studies the stability of the equilibria of \eqref{eq:coupled-oscillator} before and after bifurcation. 

\begin{proposition}[Bifurcation in oscillator networks]\label{prop:saddle-node}
Consider the nonlinear coupled oscillators~\eqref{eq:coupled-oscillator} over a weighted undirected acyclic graph $\mathcal{G}=\{\mathcal{V},\mathcal{E}, \mathcal{W}\}$. If the parameterized  natural frequencies satisfy Assumption~\ref{ass:bifurcation} for the singleton $2$-cutset $\partial S=\{e\}$, then a bifurcation occurs in~\eqref{eq:coupled-oscillator} such that
\begin{enumerate}
    \item\label{p1:subcritical} for $\alpha\in [0,\alpha^*)$, the dynamical system~\eqref{eq:coupled-oscillator} has a locally asymptotically stable equilibrium manifold $\overline{\theta}_1(\alpha)+\mathrm{span}\{\vect{1}_n\}$ given by 
    \begin{align*}
        \overline{\theta}_1(\alpha)=L^{\dagger}B\mathcal{A}\sin^{-1}(B^{\top}L^{\dagger}\omega(\alpha)),
    \end{align*}
    and an unstable equilibrium manifold $L^{\dagger}B\mathcal{A}\mathbf{n}+\mathrm{span}\{\vect{1}_n\}$ where $\mathbf{n}\in \real^{|\mathcal{E}|}$ is 
    \begin{align*}
        \mathbf{n}_s =\begin{cases}
        \sin^{-1}(B^{\top}L^{\dagger}\omega(\alpha))_s & s\ne e\\
        \pi - \sin^{-1}(B^{\top}L^{\dagger}\omega(\alpha))_s & s=e.  
\end{cases}
    \end{align*}
    \item\label{p3:supercritical} for $\alpha \in (\alpha^*, \infty)$, the dynamical system~\eqref{eq:coupled-oscillator} has no equilibrium manifold. 
\end{enumerate}
\end{proposition}

\begin{remark}
While Proposition~\ref{prop:saddle-node}\ref{p1:subcritical} implies that for $\alpha \in [0,\alpha^*)$ there exist a locally stable and a locally unstable equilibrium manifold for the coupled oscillators network~\eqref{eq:coupled-oscillator}, it does not exclude the existence of other equilibrium manifolds. 

\end{remark}

\subsection{Attraction-Repulsion Dynamics}
Consider the attraction-repulsion dynamics~\eqref{eqn:nougf86} over a network $\mathcal{G}=\{\mathcal{V},\mathcal{E}, \mathcal{W}\}$ and suppose that the attraction coefficients $w^a=\{w_e^a\}$ are parameterized by a real-valued parameter $\alpha \in [0,\infty)$ satisfying the following assumption. 

\begin{assumption}\label{ass:bifurcation_attraction}
There exist a $2$-cutset $\partial S$ for some $S\subseteq \mathcal{V}$ and an $\alpha^*\in (0, \infty)$ such that 
\begin{enumerate}
    \item If $\alpha \in (0, \alpha^*)$ then $w^a_{e}(\alpha)>w^r_{e}$, for every $e\in \mathcal{E}$;
    \item If $\alpha=\alpha^*$, then $w^a_{e}(\alpha^*)=w^r_{e}$ for all $e\in \partial S$; 
     \item If $\alpha \in (\alpha^*, \infty)$ we have
     \begin{align*}
         w^a_{e}(\alpha)&<w^r_{e}, \mbox{ for every } e\in \partial S, \\
         w^a_{e}(\alpha)&>w^r_{e}, \mbox{ for every }   e\in \mathcal{E} \setminus  \partial S.
     \end{align*}
\end{enumerate}
\end{assumption}
If Assumption \ref{ass:bifurcation_attraction} holds, then at $\alpha=\alpha^*$, the attraction and repulsion coefficients on the  $2$-cutset $\partial S$ are equal. The next proposition shows that, under this assumption, a bifurcation happens at $\alpha=\alpha^*$ which leads to instability of the consensus. 

\begin{proposition}[Bifurcation in Attraction-Repulsion Dynamics]\label{prop:attraction}
Consider the attraction-repulsion dynamics \eqref{eqn:nougf86} over a  weighted undirected graph $\mathcal{G}=\{\mathcal{V},\mathcal{E}, \mathcal{W}\}$. Then if the parameterized  attraction function satisfies Assumption~\ref{ass:bifurcation_attraction} for a 2-cutset $\partial S$ for some $S\subseteq \mathcal{V}$,
\begin{enumerate}
    \item\label{p1:ar} for $\alpha\in [0,\alpha^*)$, the  system~\eqref{eqn:nougf86} converges globally to $\frac{1}{n}\vect{1}_n\vect{1}_n^{\top}x(0)$; 
    \item\label{p2:ar} for $\alpha \in (\alpha^*,\infty)$, the invariant manifold $\mathcal{M}=\mathrm{span}\{\vect{1}_n\}$ is unstable for the system~\eqref{eqn:nougf86}. 
\end{enumerate}
\end{proposition}

\begin{remark}
Proposition~\ref{prop:attraction} does not exclude the emergence of other contracting attractors after bifurcation. Indeed, in \cite{gazii}, it is shown that for attraction-repulsion dynamics~\eqref{eqn:nougf86} states of the agents converge to a positive invariant set defined as 
\begin{align*}
    B_{\epsilon}=\{x: \|x-\overline{x}\|\leq \beta\}, \quad {\rm where} \quad \beta=\frac{w_r}{w_a}\sqrt{\frac{c}{2}}{\rm exp}(\frac{-1}{2}).
\end{align*}
\end{remark}

\section{Detection and Localization of critical transitions}

Now, we extend the deterministic and the stochastic detection algorithms in Section~\ref{sec:detection-2} to nonlinear networks.


\subsubsection*{Coupled Oscillator Networks} Let $\overline{\theta}_1(\alpha)+ \mathrm{span}\{\vect{1}_n\}$ be the locally asymptotically stable equilibrium manifold of \eqref{eq:coupled-oscillator} as described in Proposition~\ref{prop:saddle-node}. Using the change of coordinate $\epsilon = \theta-\overline{\theta}_1(\alpha)$, we can approximate the perturbation dynamics of  coupled oscillators as follows
\begin{align}\label{variational-co}
    \frac{d\epsilon}{dt} = -B\mathcal{A}_{\mathrm{CO}}(\alpha)B^{\top}\epsilon,
\end{align}
where $\mathcal{A}_{\mathrm{CO}}(\alpha) = \mathcal{A}[\cos(B^{\top}\overline{\theta}_1(\alpha))]$. 
We define the family of graphs $\mathcal{G}(\alpha)=\{\mathcal{V},\mathcal{E},\mathcal{W}(\alpha)\}$ such that 
\begin{align*}
   \left(\mathcal{W}(\alpha)\right)_{ij} =  \left(\mathcal{A}_{\mathrm{CO}}(\alpha)\right)_{ee},\qquad\mbox{ for all }\alpha \in [0,\infty),
\end{align*}
where $e=(i,j) \in \mathcal{E}$. Using Proposition~\ref{prop:saddle-node}, one can show that $B\mathcal{A}[\cos(B^{\top}\overline{\theta}_1(\alpha))]B^{\top}$ is the Laplacian matrix of the graph $\mathcal{G}(\alpha)$ and thus it is a positive semi-definite matrix. We denote the eigenvalues of this Laplaican matrix by $0=\lambda_1(\alpha)\le\lambda_2(\alpha)\le\ldots\le \lambda_n(\alpha)$ and their associated normalized eigenvectors by $\frac{1}{\sqrt{n}}\vect{1}_n=v_1(\alpha),v_2(\alpha), \ldots,v_n(\alpha)$.




\subsubsection*{Attraction-Repulsion Networks}
 Using the change of coordinate $\epsilon = x - \overline{x}$, where $\overline{x}=\frac{1}{n}\vect{1}_n\vect{1}_n^{\top} x(0)$, we can approximate the perturbation dynamics as follows
\begin{align}\label{eq:variational-AR}
    \frac{d\epsilon}{dt} = -\frac{\partial}{\partial x}\big(B \overline{\mathcal{A}}(x,\alpha)B^{\top}x \big)\bigg|_{\overline{x}}\epsilon.
\end{align}
For every $x\in \real^n$ and every $\alpha\in [0,\infty)$, we compute
\begin{align}\label{eqn:n98tg87w}
\frac{\partial}{\partial x}\big(B\overline{\mathcal{A}}(x,\alpha)B^{\top}x \big)&=- \frac{\partial}{\partial x}\big(B{\mathcal{A}_{AR}}(x,\alpha) \big)B^{\top}x \nonumber\\
&- B\overline{\mathcal{A}}(\overline{x},\alpha)B^{\top}.
\end{align}
For every $\alpha\in [0,\infty)$, the first term in the RHS of \eqref{eqn:n98tg87w} is always zero when $x=\overline{x}$ because $B^{\top}\overline{x}=0$. However, using the fact that $\overline{\mathcal{A}}(\overline{x},\alpha^*)=\vect{0}_{n\times n}$, one can show the second term in the RHS of \eqref{eqn:n98tg87w} is zero at the equilibrium point $x=\overline{x}$ only at the bifurcation,  $\alpha=\alpha^*$. Therefore, the perturbation dynamics~\eqref{eq:variational-AR} can be written as follows: 
\begin{align}\label{eq:error-AR}
    \frac{d\epsilon}{dt} = - B\mathcal{A}_{\mathrm{AR}}(\alpha)B^{\top}\epsilon.
\end{align}
where $\mathcal{A}_{\mathrm{AR}}(\alpha) = \overline{\mathcal{A}}(\overline{x},\alpha)$. We define the family of graphs $\mathcal{G}(\alpha)=\{\mathcal{V},\mathcal{E},\mathcal{W}(\alpha)\}$ such that 
\begin{align*}
   (\mathcal{W}(\alpha))_{ij} =  (\mathcal{A}_{\mathrm{AR}}(\alpha))_{ee},\qquad\mbox{ for all }\alpha \in [0,\infty),
\end{align*}
where $e=(i,j)\in \mathcal{E}$. Using Proposition~\ref{prop:attraction}, one can show that  $B\overline{\mathcal{A}}(\overline{x},\alpha)B^{\top}$ is the Laplacian matrix of the graph $\mathcal{G}(\alpha)$ and thus it is a positive semi-definite matrix. For the simplicity of notations, similar to the case of nonlinear oscillators, we denote the eigenvalues of this Laplacian matrix by $0=\lambda_1(\alpha)\le\lambda_2(\alpha)\le\ldots\le \lambda_n(\alpha)$ and the their associated normalized eigenvectors by $\frac{1}{\sqrt{n}}\vect{1}_n=v_1(\alpha),v_2(\alpha), \ldots,v_n(\alpha)$.

The spectral properties of the Laplacian matrix $B\mathcal{A}_{\mathrm{CO}}(\alpha^*)B^{\top}$ (resp. $B\mathcal{A}_{\mathrm{AR}}(\alpha^*)B^{\top}$) plays a critical role in our detection algorithm. Before we state the next lemma, recall that $\underline{\mathcal{G}}$ is the lower-bound graph associated with the family of parameterized graphs $\mathcal{G}(\alpha)$ for $\alpha\in [0,\infty)$ defined in Section~\ref{sec:definitions}. The following lemma plays a crucial role in our detection and localization algorithms. 


\begin{lemma}\label{lem:1}
Consider the perturbed coupled oscillator dynamics~\eqref{variational-co} (resp. perturbed attraction-repulsion dynamics~\eqref{eq:error-AR}) with bifurcation parameter $\alpha\in [0,\infty)$ satisfying Assumption~\ref{ass:bifurcation} (resp. Assumption~\ref{ass:bifurcation_attraction}) with a 2-cutset $\partial S$ for some $S\subseteq \mathcal{V}$. Then the following statements hold:
\begin{enumerate}
    \item\label{p1:eigenvalue} $\lim_{\alpha\to\alpha^*}\lambda_2(\alpha) = \lambda_2(\alpha^*)=0$;
    \item\label{p2:eigenvector} $\lim_{\alpha\to \alpha^*} v_2(\alpha) = v_2(\alpha^*) = \frac{1}{\sqrt{n}}\chi^S$, 
\end{enumerate}
where $\chi^S$ is the indicator vector defined in Section \ref{sec:definitions}
\end{lemma}

\subsection{Network critical slowing down: Deterministic approach}

Our key observation is that the system's recovery rate converges to zero as the it approaches the bifurcation. Moreover, the structure of the persisting perturbation determines the location of the bifurcation in the network. Now we can state the main result of this paper. 

\begin{theorem}[Convergence of perturbation dynamics]\label{thm:asymbehavior-kuramoto}
Consider the nonlinear coupled oscillator network~\eqref{eq:coupled-oscillator} (resp. attraction-repulsion network~\eqref{eqn:nougf86}) satisfying Assumption~\ref{ass:bifurcation} on an acyclic graph (resp. satisfying  Assumption~\ref{ass:bifurcation_attraction} on an arbitrary graph). For every $\alpha\in [0,\alpha^*)$ and every solution of $\epsilon:\real_{\ge 0}\to \real^n$ of the perturbation dynamics~\eqref{variational-co} (resp. perturbation dynamics~\eqref{eq:error-AR}), the following statements hold
\begin{enumerate}
    \item\label{p1:kuramoto} we have
    \begin{align*}
      \lim_{t\to \infty} \epsilon(t) = \tfrac{1}{n}\vect{1}_n\vect{1}_n^{\top}\epsilon (0)  
    \end{align*}
     with the rate of convergence $\lambda_2(\alpha)$; 
\item \label{p2:kuramoto} we have
    \begin{align*}
    \lim_{t\to \infty} \epsilon(t) = \left(\tfrac{1}{n}\vect{1}_n\vect{1}_n^{\top} + e^{-\lambda_2(\alpha)t}v_2(\alpha)v_2^{\top}(\alpha)\right)\epsilon (0)
\end{align*}
uniformly in $\alpha$ with rate of convergence $\lambda_3(\underline{\mathcal{G}})$.
\end{enumerate}
\end{theorem}

\begin{remark}[\textbf{Deterministic interpretation of Network Critical Slowing Down}]\label{rem:bn9o8tg}
 Lemma~\ref{lem:1} and Theorem~\ref{thm:asymbehavior-kuramoto}\ref{p1:kuramoto} show that the rate of convergence of $\epsilon(t)$ to the consensus decreases to zero as the system gets closer to the bifurcation, i.e., $\alpha\to \alpha^*$. This phenomenon, termed as ``network critical slowing down", is the generalization of the slowing down phenomenon discussed in Section \ref{sec:detection-2}. The uniform convergence in Theorem~\ref{thm:asymbehavior-kuramoto}\ref{p2:kuramoto} can be used to detect the bifurcation, as explained in the following subsection. 
 \end{remark}

 The proposed bifurcation detection and localization algorithm discussed in this subsection is based upon the phenomenon of critical slowing down discussed above. 
 Using the uniform convergence result in Theorem~\ref{thm:asymbehavior-kuramoto}\ref{p2:kuramoto}, for every $\alpha\in [0,\alpha^*)$ the solution $\epsilon:\real_{\ge 0}\to \real^n$ of the perturbation dynamics has the following form
    \begin{align}\label{eqn:gv9aiuwgf987}
        \epsilon(t) = \Big(\tfrac{1}{n}\vect{1}_n\vect{1}_n^{\top}+ e^{-\lambda_2(\alpha)t} v_2(\alpha)v_2^{\top}(\alpha)\Big)\epsilon (0) + \mathcal{O}(e^{-\lambda_3(\underline{\mathcal{G}}) t}).
    \end{align}
    For $\alpha$ close to $\alpha^*$ satisfying
  $ \lambda_2(\alpha) \le \tfrac{1}{2} \lambda_3(\underline{\mathcal{G}})$ and for a detection time scale $\zeta < 1$ such that $t\ge \frac{\ln(\zeta^{-2})}{\lambda_3(\underline{\mathcal{G}})}$, we have 
    \begin{align}\label{eq:order-magnitude}
        e^{- \lambda_3(\underline{\mathcal{G}}) t} &\in \mathcal{O}(\zeta^2),\nonumber \\ e^{-\lambda_2(\alpha)t} v_2(\alpha)v_2^{\top}(\alpha)\epsilon (0) &\in \mathcal{O}(\zeta). 
    \end{align}
 This means that as $t\to \infty$, the term $\mathcal{O}(e^{-\lambda_3(\underline{\mathcal{G}}) t})$ in \eqref{eqn:gv9aiuwgf987}  vanishes faster than the second term, $e^{-\lambda_2(\alpha)t} v_2(\alpha)v_2^{\top}(\alpha)\epsilon (0)$. Thus, equation \eqref{eq:order-magnitude} shows that for a considerable time period, the term $e^{-\lambda_2(\alpha)t} v_2(\alpha)v_2^{\top}(\alpha)\epsilon (0)$ persists and can be used to detect the bifurcation.

 To localize the bifurcating edge, according to Theorem \ref{thm:asymbehavior-kuramoto}, we define the vectors of residual measurement $r(t):=\epsilon(t) -\tfrac{1}{n}\vect{1}_n\vect{1}_n^{\top}\epsilon (0)$ to approximate the direction of the Fiedler vector of $\mathcal{G}(\alpha)$, denoted by $v_2(\alpha)$. It is well-known that the sign pattern of the Fiedler vector can be used to detect the graph's edge cut \cite{Fiedler}. Therefore, we can use the sign structure of the residual measurement vectors, i.e., $\mathrm{sign}(\epsilon(t)-\tfrac{1}{n}\vect{1}_n\vect{1}_n^{\top}\epsilon (0))$  to localize the bifurcating edge $e$. In particular, the bifurcating edge cut is the edge set connecting nodes with positive signs to the nodes with negative sign. The procedure is explained in Algorithm~\ref{alg:1}.


\begin{algorithm}
    \caption{\textbf{Bifurcation Detection \& Localization}}
    \label{alg:active_solver}
    \begin{algorithmic}[1]
      \REQUIRE{The initial perturbation $\epsilon(0)$, the detection time-scale $\zeta<1$, and the residual threshold $\delta$.}
      \ENSURE{Detecting bifurcation and localizing bifurcating edge cut $\partial S$.}

      \smallskip

       \STATE Set $t^* = \frac{\ln(\zeta^{-2})}{\lambda_3(\underline{\mathcal{G}})}$ 

    \STATE Collect $\{\epsilon(t)\}_{t\ge t^*}$

    \STATE Compute the residual vector $r(t) := \epsilon(t)-\tfrac{1}{n}\vect{1}_n\vect{1}_n^{\top}\epsilon(0)$ 
\IF{$\|r(t)\|_2\ge \delta$}
 \STATE An edge cut undergoes a bifurcation.\\
$S:=\{j\in \mathcal{V}: r_j(t)>0\}$\\
  $\partial S$ are the bifurcating edges.  

 \ENDIF

    \end{algorithmic}
\label{alg:1}
  \end{algorithm}
  
  \begin{remark}
  The following remarks are in order. 
  \begin{enumerate}
\item To identify the bifurcating edges based on Algorithm \ref{alg:1}, we need to have a complete knowledge of the edge cut set $\partial S$. For that, a prior knowledge about the binary network topology is required to localize the bifurcating edges. 

 \item Note that the values of the detection time-scale ($\zeta<1$), the residual threshold ($\delta$), as well as the magnitude and structure of the initial perturbation $\epsilon (0)$ must be  determined (or learned) based on the specification of the physical system, precision of measurements, and the application of interest.

 \item Algorithm~\ref{alg:active_solver} proposes a completely data-driven method for detection and localization of bifurcations. We can obtain the equilibrium manifold by  having the state measurements (i.e., $x$) over a long period of time. Then, $\epsilon=x-\overline{x}_{\alpha}$ in Algorithm~\ref{alg:active_solver} is used for detection and localization of bifurcating edges.

 \item In many applications, the initial perturbation is the last instance of an exogenous signal $u_{\mathrm{exg}}:[0,t_0]\to \real^n$ applied on the network dynamics. In this case, we directly use Theorem~\ref{thm:asymbehavior-kuramoto} for the initial condition $\epsilon(t_0)$ instead of $\epsilon(0)$. We use this setup for our simulations in the next example. 
  \end{enumerate}
  \end{remark}

\begin{example}\label{exp:1}
Consider a coupled oscillator network~\eqref{eq:coupled-oscillator} over an acyclic graph shown in Fig.~\ref{fig:example1}. 
\begin{figure}
              \centering
              \includegraphics[scale=0.35]{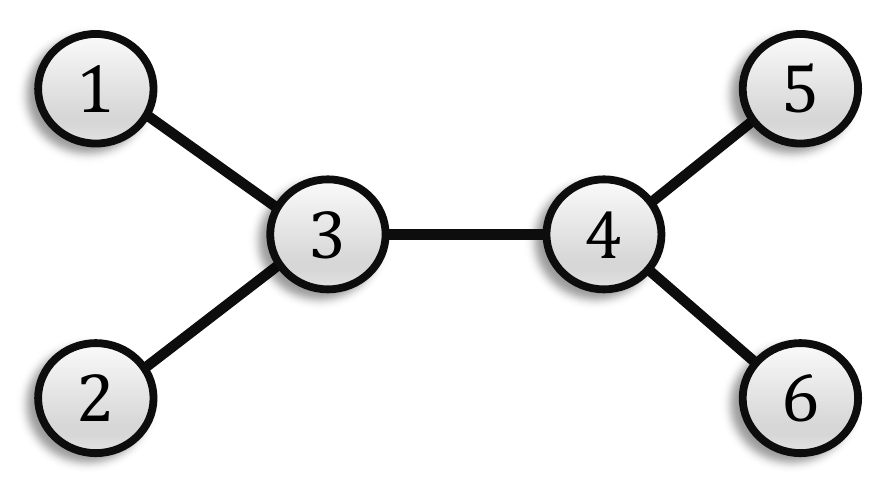}
     \caption{\footnotesize{Topology of the oscillator network in Example~\ref{exp:1}.}}
              \label{fig:example1}
\end{figure}
The objective is to show the performance of Algorithm \ref{alg:1} in detecting the bifurcation and identifying the bifurcating edge (here the edge connecting nodes 2 and 3). A perturbation signal $u_{\mathrm{exg}}=[2 \hspace{1mm} 0 \hspace{1mm} 0 \hspace{1mm} 0 \hspace{1mm} 0 \hspace{1mm} 0]^{\top}$ is applied on the network from $t=2s$ to $t=3s$. The components of the residual signal $r$, defined in Algorithm \ref{alg:1}, for a system close to the bifurcation and for a system away from the bifurcation are shown in Fig.~\ref{fig:2m16} (a) and (b), respectively. 
\begin{figure}[t!]
\centering
\includegraphics[scale=0.54]{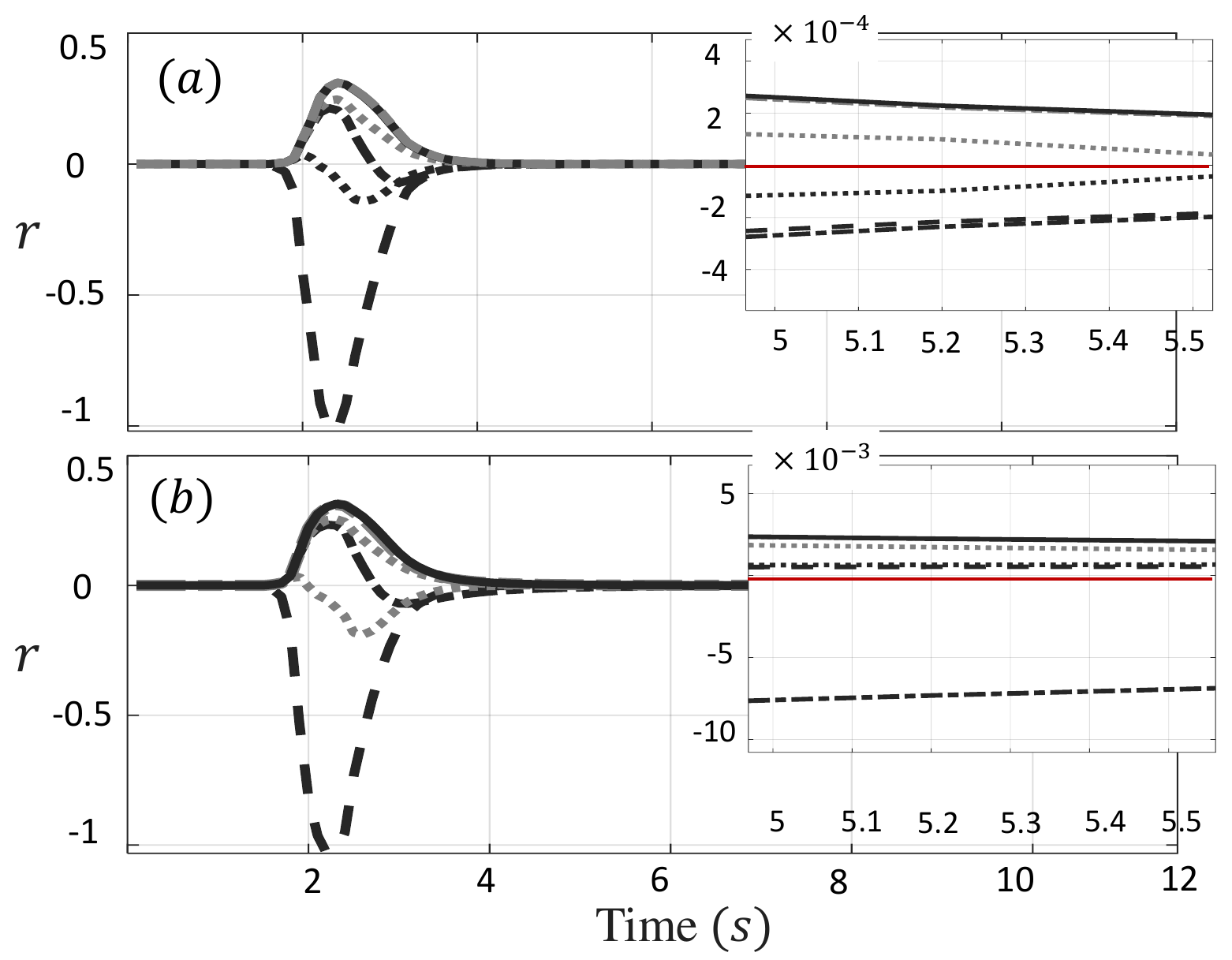}
\caption{\footnotesize{Evolution of the residual signal $r$ for coupled oscillator dynamics~\eqref{eq:coupled-oscillator} over the acyclic network~\ref{fig:example1} far from bifurcation is shown in the top plot, and close to the bifurcation is shown in the bottom plot.}}
\label{fig:2m16}
\end{figure}
By looking at the residual at $t\approx 5s$, we see that its magnitude for the system close to the bifurcation is more than ten times larger than that of the system away from bifurcation. By setting the threshold value $\delta\approx 10^{-3}$ in Algorithm \ref{alg:1} and giving enough time to the system to settle the perturbation (here $t\approx 5s$), for the system away from the bifurcation, the $2$-norm of the residual signal does not pass $\delta$ threshold and, for the system close to the bifurcation, this value passes the threshold. More precisely, we have $\|r(5)\|_2^{\rm before}=5.8\times 10^{-4}$ and $\|r(5)\|_2^{\rm after}=8.3\times 10^{-3}$. To identify the bifurcating edge, we only need to check the sign pattern of the residual measurement $r$ which shows the clustering in the network. As shown in Fig.~\ref{fig:2m16}, bottom, all elements have positive values, except the one which corresponds to node 2 in the network. This shows that the edge between nodes 2 and 3 is bifurcating. 

The algorithm is also tested on attraction-repulsion dynamics~\eqref{eqn:nougf86} and the results are shown in Fig.~\ref{fig:2m1666}. In Fig.~\ref{fig:2m1666} $(a)$, the residuals when far from bifurcation is shown and in Figs.~\ref{fig:2m1666} $(b)$ and $(c)$, the residual close to the bifurcation, when edges $(2,3)$ and $(3,4)$ are  bifurcating, respectively, are shown.
\begin{figure}[t!]
\centering
\includegraphics[scale=0.50]{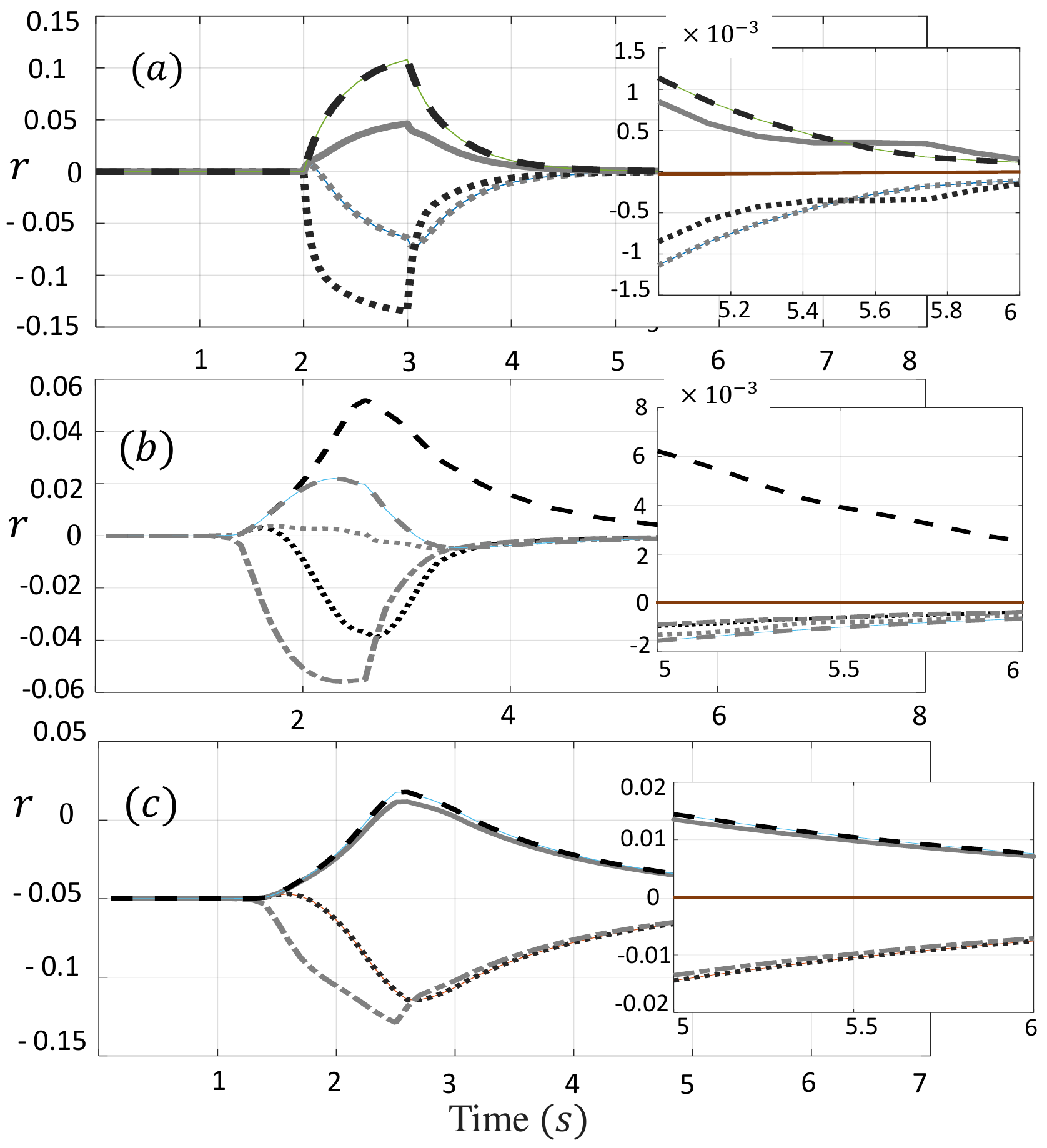}
\caption{\footnotesize{Evolution of the residual signal $r$ for attraction-repulsion dynamics~\eqref{eqn:nougf86} over the  network~\ref{fig:example1}. (a) far from bifurcation, (b) when edge $(2,3)$ becomes close to bifurcation, (c) when edge  $(3,4)$ becomes close to bifurcation.}}
\label{fig:2m1666}
\end{figure}

\end{example}

\subsection{Network critical slowing down: Stochastic approach}

We now revisit the stochastic detection method discussed in Section \ref{sec:variancee} and extend it to networks. Similar to the scalar case, we assume that stochastic disturbances are applied in specific time steps $\Delta t$, for $\Delta t\in \mathbb{Z}_{\ge 0}$ to the dynamical system. In this case, for close enough states to the $\overline{x}$, the dynamics of the unperturbed systems is approximated by \eqref{eqn:buf97g8}. The state of the system after applying the stochastic white noise $\{\xi_{t}\}_{t=1}^{\infty}$ with covariance matrix $\sigma I_n$ is given by 
\begin{align}
      x_{t+1} - \overline{x}_{\alpha} &= e^{\frac{\partial f}{\partial   x}\big|_{\overline{x}_{\alpha}}\Delta t} (  x_t-\overline{x}_{\alpha}) +   \xi_t.
\end{align}
Thus, by introducing $  e_{t} =   x_{t} -\overline{x}_{\alpha}$, we get 
\begin{align}\label{eqn:cn08g86r}
      e_{t+1} = \Gamma(\alpha)   e_{t} +   \xi_t,
\end{align}
where $\Gamma(\alpha)=e^{\frac{\partial f}{\partial x}\big|_{\overline{x}_{\alpha}}\Delta t}$. Equation~\eqref{eqn:cn08g86r} can be considered as a multivariate counterpart of autoregressive model \eqref{eqn:nc0a89v8}. We need some prepossessing on the state data before further analysis. Let $Q\in \real^{n-1\times n}$ be a matrix such that \begin{align}\label{eqn:n08g75e}
    Q Q^{\top} = I_{n-1}, \qquad Q^{\top}Q = I_n - \tfrac{1}{n}\vect{1}_n\vect{1}_n^{\top}.
\end{align}
\begin{lemma}\label{lem:v86yd6tdc}
Defining $\overline{ {e}}_{t}=Q {e}_t$ and $\overline{{\Gamma}}(\alpha)=Q \Gamma(\alpha)Q^\top$, we have 
\begin{equation}
\lim_{t\to \infty}\mathrm{tr}\big({\rm Cov}[\overline{ {e}}_t]\big)=\sum_{i=1}^{n}\frac{\sigma^2}{1-\lambda^2_i(\overline{\Gamma}(\alpha))}.
\end{equation}
\end{lemma}

The following two theorems characterize the error covariance of nonlinear coupled oscillators and attraction-repulsion dynamics at bifurcation.

\begin{theorem} \label{thm:nonlinearoscstoch}
Consider the nonlinear coupled oscillator network~\eqref{eq:coupled-oscillator} with parameterized natural frequencies satisfying Assumption~\ref{ass:bifurcation} and attraction-repulsion dynamics~\eqref{eqn:nougf86} with attraction parameters satisfying Assumption~\ref{ass:bifurcation_attraction}. Both dynamics are subjected to a stochastic white noise process $\{\xi_{t}\}_{t=1}^{\infty}$ with covariance  $\sigma I$. Then, the following statements hold:
\begin{enumerate}
    \item for $\alpha\in [0,\alpha^*)$, for both dynamics we have 
    \begin{equation}\label{eqn:before_bifurcation}
\lim_{t\to \infty}\mathrm{tr}\big({\rm Cov}[\overline{ {e}}_t]\big) < \infty.
\end{equation} 
    \item at bifurcation, i.e., $\alpha=\alpha^*$, for both dynamics we have 
\begin{equation}\label{eqn:niuf87s}
\lim_{t\to \infty}\mathrm{tr}\big({\rm Cov}[\overline{ {e}}_t]\big) = \infty.
\end{equation} 
\end{enumerate}
\end{theorem}

\begin{remark}[\textbf{Stochastic Interpretation of Network Critical Slowing Down}]\label{rem:bn9ofwef8tg}

On the verge of a bifurcation, the impacts of shocks decay very slowly and their accumulating effect increases the variance of the state variable as shown in Theorem~\ref{thm:nonlinearoscstoch}. In principle, critical slowing down could reduce the ability of the system to track the fluctuations, and thereby produce an opposite effect on the variance. 
 \end{remark}

\begin{remark}[\textbf{Decentralized Detection}]
Theorem~\ref{thm:nonlinearoscstoch} naturally leads to a stochastic algorithm for detection of bifurcations in coupled oscillator networks and attraction-repulsion dynamics. One of the key features of this algorithm is that it can be implemented in a completely decentralized fashion.  Indeed, each agent has access to its own state and can compute the corresponding variance. The trace of the covariance goes to infinity if and only if the variance of at least one node goes to infinity. Hence, if a node finds its variance going to infinity, without getting information from other nodes, it can raise an alarm for the bifurcation in the network. \label{rem:boiufgviyfdu}
\end{remark}

\section{Discussion: Comparing Deterministic and Stochastic Approaches}
 \label{sec:conclusion}
Each of the above deterministic and stochastic methods  has limitations and advantages compared to the other. 

\subsubsection{Non-Intrusive Nature and Measurement Precision}
One of the advantages of stochastic method (compared to the deterministic method) is its non-intrusive nature. Indeed, the stochastic approach to detect critical transitions relies on the additive noise, which usually exists in most physical systems.  Moreover, in the stochastic approach, the detection of the bifurcation is solely based on the growth of the state variance. Hence, even a coarse  threshold for the variance is enough to detect the critical transition. However, in the deterministic method, according to Algorithm \ref{alg:1}, some prior knowledge about the system, i.e., the values of $\delta, \zeta$, is required to detect and localize the bifurcation.

\subsubsection{Localization}
While the deterministic method can precisely localize the bifurcating edges, the stochastic approach can only detect the occurrence of a  bifurcation without localizing it in the network. In fact, this extra feature of the deterministic method in both detecting and localizing the bifurcation in networks comes with the price of being intrusive and requiring precise measurements.  


\subsubsection{Decentralized Implementation}
Unlike the deterministic approach, the stochastic detection method can be implemented in a decentralized manner. The mechanism was discussed in Remark \ref{rem:boiufgviyfdu}.

The advantages and limitations of deterministic and stochastic methods are summarized in Table \ref{tab:nscvh}. 

\begin{table}[h] 
\begin{tabular}{c c c c c} 
 \hline
\textbf{Method} & Decentralized &  Non-Intrusive & Localization \\
 \hline
Deterministic & No & No & Yes \\
Stochastic & Yes & Yes & No  \\ 
\end{tabular}
\centering
 \caption {Comparison of the two detection methods. }
\label{tab:nscvh}
\end{table}

\section{Conclusion and Future Directions}
\noindent\textbf{Summary.} This paper initiates a new research avenue in detection and identification of critical transitions in nonlinear networks. Despite the fact that nonlinear systems do not follow canonical structures, we showed that bifurcations in those systems leave signatures which are easy to detect and identify with cost-effective methods. We proposed a deterministic and a stochastic data-driven method, based on the phenomenon of critical slowing down, to detect bifurcations for some classes of nonlinear networks. Our deterministic approach is based on detecting slow recovery rate after small perturbations, and our stochastic method is based on measuring the variance of states subject to an additive process noise. The phenomenon of critical slowing down was leveraged to provide a heuristic algorithm to identify critical links under bifurcation. We compared the abilities of the two methods. 
\newline
 \textbf{Limitations and Future Work.}  Our work can be extended in several ways. Our analysis only pertained two classes of nonlinear systems. Generalization of our methods to larger classes of Laplacian-based nonlinear systems is a possible future direction. Moreover, extension to other network topologies (e.g., digraphs) is an avenue for further investigation. In this paper, we showed that near bifurcation the recovery rate tends to zero and the variance tends to infinity. Quantifying the indices  based on the distance to the bifurcation is a straightforward but an important study for the future. In the stochastic method, our analysis was focused on the second moment, i.e., variance of the states' distribution. It is shown that critical transitions also leave signatures on higher order moments, e.g., skewness and kurtosis \cite{Dakos}. Analyzing these effects in nonlinear networks is also an interesting research avenue.

\section{Acknowledgments} 
The first author thanks  Shreyas Sundaram for discussions about the phenomenon of critical slowing down and the second author thanks Francesco Bullo for discussions about nonlinear oscillator networks.

\bibliographystyle{IEEEtran}
\bibliography{Refs}

\input{Appendix}

\end{document}

%% file: Appendix.tex
\begin{appendix}

\subsection{Proof of Proposition~\ref{prop:saddle-node}}
\label{sec:appb}
\begin{proof}
For every $\alpha\in [0,\infty)$, the equilibrium points of~\eqref{eq:coupled-oscillator} are the solutions of the following algebraic equation:
\begin{align}\label{eq:algebraic-equilibrium}
     \omega(\alpha) = B\mathcal{A}\sin(B^{\top}\theta). 
\end{align}
Thus, by multiplying both side of the algebraic equation~\eqref{eq:algebraic-equilibrium} in $B^{\top}L^{\dagger}$, we get
\begin{align}\label{eq:saberissmart}
    B^{\top}L^{\dagger}\omega(\alpha) = B^{\top}L^{\dagger}B\mathcal{A}\sin(B^{\top}\theta)=\sin(B^{\top}\theta),
\end{align}
where the last equality holds because we have $B^{\top}L^{\dagger}B\mathcal{A}=I_{n}$ for acyclic graphs~\cite[Theorem 5]{saber}. Regarding part~\ref{p1:subcritical}, for $\alpha\in [0,\alpha^*)$, by Assumption~\ref{ass:bifurcation}, we get $\|B^{\top}L^{\dagger}\omega(\alpha)\|_{\infty}<1$. Therefore solution of~\eqref{eq:saberissmart} should satisfy 
\begin{align}\label{eq:saberisgood}
    (B^{\top}\theta(\alpha))_s &  = \sin^{-1}(B^{\top}L^{\dagger}\omega(\alpha))_s, \mbox{ or} \nonumber\\
    (B^{\top}\theta(\alpha))_s & = \pi -\sin^{-1}(B^{\top}L^{\dagger}\omega(\alpha))_s, 
\end{align}
for every $s\in \mathcal{E}$. Using some algebraic manipulations, one we can find two equilibrium manifolds $\theta_1(\alpha)+\mathrm{span}\{\vect{1}_n\}$ and $\theta_2(\alpha)+\mathrm{span}\{\vect{1}_n\}$ for the dynamics~\eqref{eq:coupled-oscillator} as:
\begin{align*}
    \theta_1(\alpha) &= L^{\dagger}B\mathcal{A}\sin^{-1}(B^{\top}L^{\dagger}\omega(\alpha));\\
    \theta_2(\alpha) &= L^{\dagger}B\mathcal{A} \mathbf{n}, 
\end{align*}
where $\mathbf{n}\in \real^{|\mathcal{E}|}$ is defined as follows:
    \begin{align*}
        \mathbf{n}_s =\begin{cases}
        \sin^{-1}(B^{\top}L^{\dagger}\omega(\alpha))_s & s\ne e\\
        \pi - \sin^{-1}(B^{\top}L^{\dagger}\omega(\alpha))_s & s=e.  
\end{cases}
    \end{align*}
Regarding the local stability of the manifolds $\theta_1(\alpha)+\mathrm{span}\{\vect{1}_n\}$ and $\theta_2(\alpha)+\mathrm{span}\{\vect{1}_n\}$, we compute the Jacobian of the system~\eqref{eq:coupled-oscillator}:
\begin{align*}
    J_\theta = -B \mathcal{A}[\cos(B^{\top}\theta)]B^{\top}.
\end{align*}
At $\theta=\theta_1(\alpha)$, we have $J_{\theta_1(\alpha)} = B \mathcal{A}[\cos(B^{\top}\theta_1(\alpha))]B^{\top}$. In turn, we have $\frac{-\pi}{2}\vect{1}_m< B^{\top}\theta_1(\alpha)< \frac{\pi}{2}\vect{1}_m$ and therefore $\cos(B^{\top}\theta_1(\alpha))> 0$. This implies that $[\cos(B^{\top}\theta)]$ is a positive diagonal matrix and therefore $J_{\theta_1}$ is a negative semi-definite matrix with kernel $\mathrm{span}\{\vect{1}_n\}$. This means that $\theta_1(\alpha)+\mathrm{span}\{\vect{1}_n\}$ is locally asymptotically stable equilibrium manifold of the system. 

At $\theta=\theta_2(\alpha)$, we have $J_{\theta_2(\alpha)} = B \mathcal{A}[\cos(B^{\top}\theta_2(\alpha))]B^{\top}$. Note that, for the edge $e$, we have $\frac{\pi}{2}< (B^{\top}\theta_2(\alpha))_e<\pi$. This implies that $\cos(B^{\top}\theta_2(\alpha))_e<0$. Let $\mathbf{e}_e\in \real^{|\mathcal{E}|}$ be the $e$th standard basis of $\real^{|\mathcal{E}|}$. Since the graph $G=\{\mathcal{V},\mathcal{E},\mathcal{W}\}$ is acyclic, there exists $v\in \real^n$ such that $B^{\top}v=\mathbf{e}_e$. Therefore, 
\begin{align*}
    v^{\top}J_{\theta_2(\alpha)}v &= -\mathbf{e}_e^{\top}\mathcal{A}[\cos(B^{\top}\theta_2(\alpha))]\mathbf{e}_e \\ & = -\mathbf{e}_e^{\top}\mathcal{A} \cos((B^{\top}\theta_2(\alpha))_e) >0.  
\end{align*}
As a result, $\theta_2(\alpha)+\mathrm{span}\{\vect{1}_n\}$ is an unstable equilibrium manifold of the system. 


Regarding part~\ref{p3:supercritical}, suppose that $\theta^*(\alpha)$ is an equilibrium point of the system, for $\alpha\in (\alpha^*,\infty)$. Then by~\eqref{eq:saberissmart}, we have $B^{\top}L^{\dagger}\omega(\alpha) = \sin(B^{\top}\theta^*)$. This implies that, $\|B^{\top}L^{\dagger}\omega(\alpha)\|_{\infty}\le 1$ and therefore, we get $|(B^{\top}L^{\dagger}\omega(\alpha))_e|\le 1$. However, this is in a contradiction with Assumption~\eqref{ass:bifurcation}. Thus the network~\eqref{eq:coupled-oscillator} has no equilibrium manifold. 
\end{proof}

\subsection{Proof of Proposition~\ref{prop:attraction}}
\label{sec:appb}
\begin{proof}
Regarding part~\ref{p1:ar}, for $\alpha\in [0,\alpha^*)$, we have $w_a^{ij}>w_r^{ij}$ for all $(i,j)\in \mathcal{E}$. Thus, for each diagonal element of matrix $\bar{\mathcal{A}}$, we have $w_a^{ij}(\alpha)-w_r \exp (\tfrac{-\|x_i-x_j\|^2}{c})>0.$
Using $V=x^{\top}x$ as a Lyapunov function, we get
\begin{align*}
    \dot{V}=-x^{\top}B\bar{\mathcal{A}}B^{\top}x\leq 0.
\end{align*}
Note that $\{x\in \real^n \mid \dot{V}(x)=0\}=\mathrm{span}\{\vect{1}_n\}$ is an invariant set of the system~\eqref{eqn:nougf86}. Now, using the LaSalle's invariance principle~\cite{khalil}, every trajectory of the system~\eqref{eqn:nougf86} converges to $\mathrm{span}\{\vect{1}_n\}$. Moreover, it is easy to see that $\frac{d}{dt}\vect{1}^{\top}_n x(t) = 0$ and thus $\vect{1}_n^{\top}x(t) = \vect{1}_n^{\top}x(0)$, for every $t\in \real_{\ge 0}$. This implies that the trajectory $t\mapsto x(t)$ converges to $\frac{1}{n}\vect{1}_n\vect{1}_n^{\top}x(0)$.

Regarding part \ref{p2:ar}, we show that $\mathrm{span}\{\vect{1}_n\}$ is an unstable equilibrium manifold for the attraction-repulsion dynamics~\ref{eqn:nougf86} when $\alpha\in (\alpha^*,\infty)$. We first calculate the Jacobian \eqref{eqn:n98tg87w} at the equilibrium manifold $\mathrm{span}\{\mathbb{1}_n\}$. Note that, the first term in the RHS of the equation~\eqref{eqn:n98tg87w} is always zero on the manifold $\mathrm{span}\{\mathbb{1}_n\}$. Let $x\in \real^{|\mathcal{V}|}$ be given by 
\begin{align*}
    x_i =\begin{cases}
    1 & i\in S,\\
    0 & i\ne S
    \end{cases}
\end{align*}
Then, one can easily show that $B^{\top}x\in \real^{|\mathcal{E}|}$ is given by 
\begin{align*}
    (B^{\top}x)_e =\begin{cases}
    1 & e\in \partial S,\\
    0 & e\ne \partial S
    \end{cases}
\end{align*}
As a result, for every $x\in \mathrm{span}\{\vect{1}_n\}$ and every $\alpha\in (\alpha^*,\infty)$, 
\begin{align*}
  -x^{\top}B\bar{\mathcal{A}}(x,\alpha)B^{\top}x &= \sum_{e\in \partial S} \bar{\mathcal{A}}_{ee}(x,\alpha) \\ &= \sum_{e=(i,j)} -w^a_{ij}(\alpha) + w^r_{ij} \ge 0,
\end{align*}
where the last equality holds by Assumption~\ref{ass:bifurcation_attraction}. This implies that $\mathrm{span}\{\vect{1}_n\}$ is an unstable equilibrium manifold. \end{proof}

\subsection{Proof of Lemma~\ref{lem:1}}
\begin{proof}
 We prove the result for nonlinear coupled oscillator networks. The proof for attraction-repulsion networks is similar and we omit it for the sake of brevity. Consider the perturbation dynamics For the nonlinear coupled oscillator networks~\eqref{variational-co}. At the bifurcation, i.e., $\alpha=\alpha^*$, by Assumption~\ref{ass:bifurcation},
\begin{align*}
    |(B^{\top}L^{\dagger}\omega(\alpha^*))_e|= 1,\qquad \mbox{for every } e\in \partial S
\end{align*}
Using Proposition~\ref{prop:saddle-node}, for every $e\in \partial S$, we get
\begin{align*}
    \cos((B^{\top}\overline{\theta}_1(\alpha^*))_e) = \cos(\sin^{-1}((B^{\top}L^{\dagger}\omega(\alpha^*))_e)= 0.
\end{align*}
Therefore, for every $e=(i,j)\in \partial S$, 
\begin{align}\label{eq:easyto}
    (\mathcal{A}_{\mathrm{CO}}(\alpha^*))_{ij} = \cos(B^{\top}\overline{\theta}_1(\alpha)_e) = 0.
\end{align}
This means that $\mathrm{rank}\left(\mathcal{A}_{\mathrm{CO}}(\alpha^*)\right)=n-1$ and therefore $\mathrm{rank}\left(B\mathcal{A}_{\mathrm{CO}}(\alpha^*)B^{\top}\right) =n-2$. This implies that $\lambda_2(\alpha^*)=0$. Moreover, one can easily check that
\begin{align*}
    (B^{\top}\chi^S)_e =\begin{cases}
    1 & e\in \partial S\\
    0 & e\not\in \partial S.
    \end{cases}
\end{align*}
Using~\eqref{eq:easyto} and the above equation, we get $$\tfrac{1}{\sqrt{n}}B\mathcal{A}_{\mathrm{CO}}(\alpha^*)B^{\top} \chi^S = \vect{0}_n.$$ 
This implies that $v_{2}(\alpha^*) = \frac{1}{\sqrt{n}}\chi^S$ is an eigenvector of $B\mathcal{A}_{\mathrm{CO}}(\alpha^*)B^{\top}$ associated with the eigenvalue $\lambda_2(\alpha^*)$. By~\cite[Lemma 4.3]{DS-JS:02} the eigenvalues and eigenvectors of a symmetric matrix are continuous in parameter $\alpha$ and therefore $\lim_{\alpha\to \alpha^*} \lambda_2(\alpha)=0$ and $\lim_{\alpha\to \alpha^*} v_2(\alpha^*) =\frac{1}{\sqrt{n}}\chi^S$. \end{proof}

\subsection{Proof of Theorem~\ref{thm:asymbehavior-kuramoto}}\label{sec:appb}
\begin{proof}

First note that, for coupled-oscillator network~\eqref{eq:coupled-oscillator} (resp. attraction-repulsion network~\eqref{eqn:nougf86}), for every $\alpha\in [0,\infty)$ and $t\ge 0$, the solution of the perturbation dynamics~\eqref{variational-co} (resp. the perturbation dynamics~\eqref{eq:error-AR}) is given by:
\begin{align}\label{eqn:btc75cd}
    \epsilon(t) = \sum_{i=1}^{n} e^{-\lambda_i(\alpha)t} v_i(\alpha)v_i^{\top}(\alpha) \epsilon(0)
\end{align}
where $0=\lambda_1(\alpha),\ldots,\lambda_n(\alpha)$ are the eigenvalues with the associated normalized eigenvectors $\tfrac{1}{\sqrt{n}}\vect{1}_n=v_1(\alpha),\ldots,v_n(\alpha)$ for the matrix $B\mathcal{A}[\cos(B^{\top}\theta(\alpha))]B^{\top}$ (resp. the matrix $B\overline{\mathcal{A}}(x,\alpha)B^{\top}$). 

First, we focus on the coupled-oscillator network~\eqref{eq:coupled-oscillator} and its associated perturbation dynamics~\eqref{variational-co}. Regarding part~\ref{p1:kuramoto}, let $\epsilon:\real_{\ge 0}\to \real^n$ be the solution of~\eqref{variational-co} for $\alpha\in [0,\alpha^*]$. Using Proposition~\ref{prop:saddle-node}\ref{p1:subcritical}, for every $\alpha\in [0,\alpha^*)$,
\begin{align*}
    B^{\top}\overline{\theta}_1(\alpha) = \sin^{-1}(B^{\top}L^{\dagger}\omega(\alpha))\in \left(-\tfrac{\pi}{2}\vect{1}_m,\tfrac{\pi}{2}\vect{1}_m\right).
\end{align*}
This implies that $\cos(B^{\top}\theta_1(\alpha))>\vect{0}_m$ and therefore $0=\lambda_1(\alpha)<\lambda_2(\alpha)\le \ldots\le \lambda_n(\alpha)$, for every $\alpha\in [0,\alpha^*)$. As a result, using the formula~\eqref{eqn:btc75cd}, we get $\lim_{t\to\infty} \epsilon(t) =\tfrac{1}{n} \vect{1}_n\vect{1}_n^{\top}\epsilon(0)$.
Regarding part~\ref{p2:kuramoto}, using the formula~\eqref{eqn:btc75cd}, for every $t\ge 0$,
\begin{multline}\label{eqn:boug97f86d}
    \epsilon(t) = \tfrac{1}{n}\vect{1}_n\vect{1}_n^{\top}\epsilon(0) + e^{-\lambda_2(\alpha)t} v_2(\alpha)v^{\top}_2(\alpha)\epsilon(0) \\ + \sum_{i=3}^{n} e^{-\lambda_i(\alpha)t} v_i(\alpha)v_i^{\top}(\alpha) \epsilon(0).
\end{multline}
Since the singleton $\{e\}$ in a 2-cutset and using Lemma~\ref{lem:1}, we can deduce that $0<\lambda_3(\alpha)$, for every $\alpha\in [0,\alpha^*]$. By compactness of the set $[0,\alpha^*]$ and continuity of $\lambda_3(\alpha)$ on the parameter $\alpha$~\cite[Lemma 4.3]{DS-JS:02}, we get
\begin{align*}
    0<\lambda_3(\underline{\mathcal{G}}) \le \lambda_3(\alpha)\le \ldots\le \lambda_n(\alpha),\quad\mbox{ for all }\alpha\in [0,\alpha^*].
\end{align*}
Therefore, for every $t\ge 0$, we get
\begin{multline*}
    \big\|\epsilon(t) - \tfrac{1}{n}\vect{1}_n\vect{1}_n^{\top}\epsilon(0) - e^{-\lambda_2(\alpha)t} v_2(\alpha)v^{\top}_2(\alpha)\epsilon(0)\big\|_2 \\  \le e^{-\lambda_3(\underline{\mathcal{G}})t} \sum_{i=3}^{n}  \|v_i(\alpha)v_i^{\top}(\alpha)\|_2\|\epsilon(0)\|_2.
\end{multline*}
Note that $v_i(\alpha)$ is normalized for every $\alpha\in [0,\alpha^*]$. This implies that $\|v_i(\alpha)v_i^{\top}(\alpha)\|_2 = 1$. As a result, we get 
\begin{multline*}
    \big\|\epsilon(t) - \tfrac{1}{n}\vect{1}_n\vect{1}_n^{\top}\epsilon(0) - e^{-\lambda_2(\alpha)t} v_2(\alpha)v^{\top}_2(\alpha)\epsilon(0)\big\|_2 \\ \le  n e^{-\lambda_3(\underline{\mathcal{G}})t}\|\epsilon(0)\|_2. 
\end{multline*}
Thus it is clear that as we take $t\to \infty$, $\epsilon(t)$ converges to  $\tfrac{1}{n}\vect{1}_n\vect{1}_n^{\top}\epsilon(0) + e^{-\lambda_2(\alpha)t} v_2(\alpha)v^{\top}_2(\alpha)\epsilon(0)$ uniformly in $\alpha$ with exponential convergence rate $\lambda_3(\underline{\mathcal{G}})$. 

 
Second, we focus on attraction-repulsion networks~\eqref{eqn:nougf86} and the associated perturbed dynamics~\eqref{eq:error-AR}. Regarding part \ref{p1:kuramoto}, using Proposition \ref{prop:attraction}, we have $\lim_{t\to \infty} \epsilon(t) = \tfrac{1}{n}\vect{1}_n\vect{1}_n^{\top}\epsilon(0)$. 
For part \ref{p2:kuramoto}, the proof is similar to the coupled oscillator network case and we omit it for the sake of brevity.  
\end{proof}

\subsection{Proof of Lemma~\ref{lem:v86yd6tdc}} \label{sec:appc}
We prove the result for 
We rewrite  \eqref{eqn:cn08g86r} as
\begin{align*}
\overline{ {e}}_{t+1}=Q {\Gamma}(\alpha)Q^{\top}\overline{ {e}}_{t}+Q {\xi}_t = \overline{\Gamma}(\alpha)\overline{e}_t + Q\xi_t
\end{align*}
The covariance of state $\overline{ {e}}_t$ becomes
\begin{align}\label{eqn:n97f76458}
{\rm Cov}[\overline{ {e}}_t]&={\rm Cov}[\overline{ {\Gamma}}^{t}\overline{ {e}}_0+\sum_{j=0}^{t-1}\overline{ {\Gamma}}^{j}Q\xi_{t-j-1}]\nonumber\\
&=\overline{ {\Gamma}}^{t}{\rm Cov}[\overline{ {e}}_0]\overline{ {\Gamma}}^{t}+ \sum_{j=0}^t\overline{ {\Gamma}}^{j-1}Q.{\rm Cov}[\xi_{t-j-1}].Q^{\top}\overline{ {\Gamma}}^{j-1}\nonumber\\
&=\overline{ {\Gamma}}^{t}{\rm Cov}[\overline{ {e}}_0]\overline{ {\Gamma}}^{t}+\sigma^2\sum_{j=0}^t\overline{ {\Gamma}}^{j-1}QQ^{\top}\overline{ {\Gamma}}^{j-1}\nonumber\\ & =\overline{ {\Gamma}}^{t}{\rm Cov}[\overline{ {e}}_0]\overline{ {\Gamma}}^{t}+ \sigma^2\sum_{j=0}^t\overline{ {\Gamma}}^{2(j-1)},
\end{align} 
In the last line of \eqref{eqn:n97f76458}, we used the fact that $QQ^{\top}=I_{N-1}$ and the $\overline{ {\Gamma}}$ is symmetric. Therefore, in the time limit ($t\to \infty$) for $\alpha < \alpha^*$, it becomes
$ \lim_{t\to \infty}{\rm Cov}[\overline{ {e}}_t]=\sigma^2\big(I-\overline{ {\Gamma}}^{2} \big)^{-1}$ which gives the result.

\subsection{Proof of Theorem~\ref{thm:nonlinearoscstoch}}
\label{sec:appd}
Regarding the nonlinear coupled oscillator network~\eqref{eq:coupled-oscillator}. Considering $\frac{\partial f}{\partial x}\bigg|_{\overline{x}_{\alpha}} = -B\mathcal{A}[\cos(B^{\top}\overline{\theta}_1(\alpha))]B^{\top}$,   we can write 
\begin{align*}
  {\Gamma}(\alpha) &= \tfrac{1}{n}\vect{1}_n\vect{1}_n^{\top} + \sum_{i=2}^{n} e^{-\lambda_i(\alpha)\Delta t} v_i(\alpha)v_i^{\top}(\alpha)\\ 
\overline{ {\Gamma}}(\alpha) &= \sum_{i=2}^{n} e^{-\lambda_i(\alpha)\Delta t}Q v_i(\alpha) v^{\top}_i(\alpha) Q^{\top}  
\end{align*}
Based on Lemma~\ref{lem:1}, we have $\lim_{\alpha\to\alpha^*}\lambda_2(\alpha)=\lambda_2(\alpha^*)=0$. Moreover, as $t\to \infty$, the matrix $\overline{ {\Gamma}}(\alpha^*)$ is given by 
\begin{align*}
\lim_{t\to\infty}\overline{ {\Gamma}}(\alpha^*)= Qv_2(\alpha^*)v_2^{\top}(\alpha^*)Q^{\top}   
\end{align*}
 Using this equation along with \eqref{eqn:n08g75e}, we can verify that $\overline{ {\Gamma}}(\alpha^*)$ has an eigenvalue $\lambda_0=1$ corresponding to  eigenvector $Qv_2$. Therefore, according to Lemma \ref{lem:v86yd6tdc}, at $\alpha=\alpha^*$ we have
\begin{equation*}
\lim_{t\to \infty}\mathrm{tr}\big({\rm Cov}[\overline{ {e}}_t]\big)=\sum_{i=1}^{n}\frac{\sigma^2}{1-\lambda_i^2(\overline{ {\Gamma}}(\alpha^*))}=\infty.
\end{equation*}

For attraction-repulsion dynamics~\eqref{eqn:nougf86} with the equilibrium manifold $\overline{x}_{\alpha}=\mathrm{span}\{\vect{1}_n\}$, we have  $\frac{\partial f}{\partial x}\bigg|_{\overline{x}_{\alpha}} = (w_r-w_a)B B^{\top}$ which gives
\begin{align*}
   {\Gamma}(\alpha) &= \tfrac{1}{n}\vect{1}_n\vect{1}_n^{\top} + \sum_{i=2}^{n} e^{(n-1)(w_r-w_a)t}v_i(\alpha) v^{\top}_i(\alpha)\\ 
\overline{ {\Gamma}}(\alpha) &= \sum_{i=2}^{n} e^{(n-1)(w_r-w_a)t}Q v_i(\alpha) v^{\top}_i(\alpha) Q^{\top}  
\end{align*}
When bifurcation occurs, i.e., $\alpha=\alpha^*$, we get 
\begin{align*}
\overline{ {\Gamma}}(\alpha^*) = \sum_{i=2}^{n}Q v_i(\alpha^*) v^{\top}_i(\alpha^*) Q^{\top}    
\end{align*}
It is easy to show that, for every $i\in \{2,\ldots,N\}$,  $Qv_i$ is an eigenvector of $\overline{ {\Gamma}}$ associated with eigenvalue $1$. Therefore, by Lemma~\ref{lem:v86yd6tdc}, we have
\begin{equation*}
\lim_{t\to \infty}\mathrm{tr}\big({\rm Cov}[\overline{ {e}}_t(\alpha^*)]\big)=\sum_{i=1}^{n}\frac{\sigma^2}{1-\lambda_i^2(\overline{ {\Gamma}}(\alpha))}=\infty.
\end{equation*}

\end{appendix}